\documentclass[reqno]{amsart}

\usepackage{amsmath}
\usepackage{amssymb}
\usepackage{amsfonts}
\usepackage{amsthm}

\usepackage[utf8]{inputenc}
\usepackage[T1]{fontenc}

\usepackage{dsfont}

\usepackage[%
cal=cm,
frak=euler
]
{mathalfa}

\usepackage[dvipsnames,svgnames]{xcolor}
\colorlet{MyBlue}{DodgerBlue!50!Black}

\usepackage{acronym}
\usepackage{subfigure}

\usepackage{acronym}
\usepackage{latexsym}
\usepackage{paralist}
\usepackage{xspace}

\usepackage[numbers,sort&compress]{natbib}

\bibpunct[, ]{[}{]}{,}{n}{}{,}
\newcommand{\citeor}[2][]{\citeauthor{#2} \cite[#1]{#2}}
\newcommand{\citeorp}[2][]{\textup(\citeor[#1]{#2}\textup)}

\usepackage{hyperref}
\hypersetup{
colorlinks=true,
linktocpage=true,
pdfstartview=FitH,
breaklinks=true,
pdfpagemode=UseNone,
pageanchor=true,
pdfpagemode=UseOutlines,
plainpages=false,
bookmarksnumbered,
bookmarksopen=false,
bookmarksopenlevel=1,
hypertexnames=true,
pdfhighlight=/O,
urlcolor=Maroon,linkcolor=MyBlue!70!black,citecolor=DarkGreen!80!black, 
pdftitle={},
pdfauthor={},
pdfsubject={},
pdfkeywords={},
pdfcreator={pdfLaTeX},
pdfproducer={LaTeX with hyperref}
}

\newcommand{\R}{\mathbb{R}}

\DeclareMathOperator{\argdot}{\cdotp}
\DeclareMathOperator*{\argmax}{arg\,max}
\DeclareMathOperator*{\argmin}{arg\,min}
\DeclareMathOperator{\bd}{bd}
\DeclareMathOperator{\bigoh}{\mathcal O}

\DeclareMathOperator{\cl}{cl}

\DeclareMathOperator{\dom}{dom}

\DeclareMathOperator{\exclude}{\backslash}

\DeclareMathOperator{\hess}{Hess}
\DeclareMathOperator{\im}{im}

\DeclareMathOperator{\inter}{int}

\DeclareMathOperator{\proj}{proj}

\DeclareMathOperator{\relint}{rel\,int}

\DeclareMathOperator{\supp}{supp}

\newcommand{\dd}{\:d}

\newcommand{\eps}{\varepsilon}
\newcommand{\from}{\colon}

\newcommand{\pd}{\partial}
\newcommand{\simplex}{\Delta}

\newcommand{\abs}[1]{\left\lvert #1 \right\rvert}
\newcommand{\norm}[1]{\left\| #1 \right\|}

\newcommand{\braket}[2]{\left\langle #1 \middle\vert  #2 \right\rangle}

\newcommand{\defeq}{\equiv}
\newcommand{\eqdef}{\equiv}

\newcommand{\dis}{\displaystyle}
\newcommand{\txs}{\textstyle}

\newcommand{\insum}{\sum\nolimits}
\newcommand{\inprod}{\prod\nolimits}

\usepackage[textwidth=30mm]{todonotes}

\newcommand{\revised}[1]{#1}

\theoremstyle{plain}
\newtheorem{theorem}{Theorem}
\newtheorem{corollary}[theorem]{Corollary}
\newtheorem*{corollary*}{Corollary}
\newtheorem{lemma}[theorem]{Lemma}
\newtheorem{proposition}[theorem]{Proposition}

\theoremstyle{definition}
\newtheorem{definition}[theorem]{Definition}
\newtheorem*{definition*}{Definition}

\theoremstyle{remark}
\newtheorem{remark}{Remark}
\newtheorem*{remark*}{Remark}
\newtheorem{example}{Example}

\newenvironment{Proof}[1][Proof]{\begin{proof}[#1]}{\end{proof}}

\numberwithin{equation}{section}
\numberwithin{theorem}{section}
\numberwithin{remark}{section}
\numberwithin{example}{section}

\newcommand{\play}{\mathcal{N}}
\newcommand{\act}{\mathcal{A}}
\newcommand{\pay}{u}

\newcommand{\strat}{\mathcal{X}}
\newcommand{\game}{\mathcal{G}}
\newcommand{\brep}{\choice^{0}}
\newcommand{\BREP}{\mathtt{br}}
\newcommand{\eq}{x^{\ast}}

\newcommand{\subd}{\pd}
\newcommand{\gibbs}{G}
\newcommand{\eucl}{\Pi}

\newcommand{\primal}{\mathcal{V}}
\newcommand{\dual}{\primal^{\ast}}

\newcommand{\set}{\mathcal{S}}
\newcommand{\interval}{I}

\newcommand{\bvec}{e}

\newcommand{\choice}{Q}
\newcommand{\fench}{F}
\newcommand{\breg}{D}
\newcommand{\dkl}{\breg_{\textup{KL}}}

\newcommand{\intstrat}{\strat^{\circ}}
\newcommand{\intsimplex}{\simplex^{\!\circ}}

\newcommand{\payv}{v}

\newcommand{\score}{y}
\newcommand{\Score}{Y}

\newcommand{\temp}{\eta}
\newcommand{\tempinv}{\gamma}
\newcommand{\rate}{\phi}

\newcommand{\corr}{\chi}
\newcommand{\widebar}{\bar}

\begin{document}


\title
[Learning in games via reinforcement and regularization]
{Learning in games via\\
reinforcement and regularization}

\author[P.~Mertikopoulos]{Panayotis Mertikopoulos}
\address
[Panayotis Mertikopoulos]
{CNRS (French National Center for Scientific Research), LIG, F-38000 Grenoble, France\\
and
Univ. Grenoble Alpes, LIG, F-38000 Grenoble, France}
\email{\href{mailto:panayotis.mertikopoulos@imag.fr}{panayotis.mertikopoulos@imag.fr}}
\urladdr{\url{http://mescal.imag.fr/membres/panayotis.mertikopoulos}}

\author[W.~H.~Sandholm]{William H. Sandholm}
\address
[William H.~Sandholm]
{Department of Economics, University of Wisconsin, 1180 Observatory Drive, Madison WI 53706, USA}
\email{\href{mailto:whs@ssc.wisc.edu}{whs@ssc.wisc.edu}}
\urladdr{\url{http://www.ssc.wisc.edu/~whs/}}

\thanks{The authors are grateful to Josef Hofbauer, Joon Kwon, Rida Laraki, two anonymous referees, and seminar audiences at the Hausdorff Research Institute for Mathematics and the University of Wisconsin for many interesting discussions.}

\thanks{%
Part of this work was carried out during the authors' visit to the Hausdorff Research Institute for Mathematics at the University of Bonn in the framework of the Trimester Program ``Stochastic Dynamics in Economics and Finance''.
PM is grateful for financial support from the French National Research Agency under grant ANR\textendash GAGA\textendash13\textendash JS01\textendash 0004\textendash 01,
and the French National Center for Scientific Research under grant PEPS\textendash GATHERING\textendash 2014.
WHS is grateful for financial support under NSF Grant SES\textendash1155135.}


\newacro{ESS}{evolutionarily stable state}
\newacro{HR}{Hess\-i\-an Rie\-man\-ni\-an}
\newacro{QR}{quantal response}
\newacro{KKT}{Ka\-rush\textendash Kuhn\textendash Tuc\-ker}
\newacro{RPS}[\textsf{RPS}]{Rock-Paper-Scissors}
\newacro{MP}{Matching Pennies}
\newacro{ODE}{ordinary differential equation}
\newacro{OD}[O/D]{origin\textendash destination}
\newacro{OMD}{online mirror descent}
\newacro{FTRL}{follow the regularized leader}
\newacro{MD}{mirror descent}
\newacro{KL}{Kullback\textendash Leibler}
\newacro{URL}{unpenalized reinforcement learning}
\newacro{CFP}{continuous-time fictitious play}

\begin{abstract}
%
%
We investigate a class of reinforcement learning dynamics where players adjust their strategies based on their actions' cumulative payoffs over time \textendash\ specifically, by playing mixed strategies that maximize their expected cumulative payoff minus a regularization term.
A widely studied example is exponential reinforcement learning, a process induced by an entropic regularization term which leads mixed strategies to evolve according to the replicator dynamics.
However, in contrast to the class of regularization functions used to define smooth best responses in models of stochastic fictitious play, the functions used in this paper need not be infinitely steep at the boundary of the simplex;
in fact, dropping this requirement gives rise to an important dichotomy between steep and nonsteep cases.
In this general framework, we extend several properties of exponential learning, including the elimination of dominated strategies, the asymptotic stability of strict Nash equilibria, and the convergence of time-averaged trajectories
\revised{in zero-sum games with an interior Nash equilibrium.}
\end{abstract}

\maketitle

\setcounter{tocdepth}{1}
\vspace{-1.5em}
\tableofcontents


\section{Introduction}
\label{sec:introduction}

``Reinforcement learning'' has become a catch-all term for learning in recurring decision processes where the agents' future choice probabilities are shaped by information about past payoffs.
In game theory and online optimization, most work under this name has focused on multi-armed bandit problems,
games against Nature (adversarial or otherwise),
or simultaneous decision-making by strategically interacting players
\textendash\ for a panoramic introduction, see \citeor{SB98}.
Accordingly, reinforcement learning in games typically revolves around discrete-time stochastic processes with the stochasticity arising at least in part from the agents' randomized choices \textendash\ see e.g. \citeor{BS97}, \citeor{ER98}, \citeor{FL98}, \citeor{FS99}, \citeor{Hop99a}, \citeor{HMC00}, \citeor{Beg05}, \citeor{LC05}, \citeor{CMS10}, \citeor{CGM15} and many others.

A key approach to analyzing these processes is the \ac{ODE} method of stochastic approximation, a method which relates the behavior of the stochastic model under study to that of a ``mean field'' \ac{ODE} \citeorp{Ben99}.
Motivated by the success of this approach, we follow \citeor{Sor09} and \citeor{HSV09} by specifying a reinforcement learning scheme directly in continuous time.
In so doing, we are able to focus squarely on the deep relations between reinforcement learning, convex analysis and population dynamics.%
\footnote{One could use the results developed in this paper to analyze discrete-time learning schemes as in \citeor{LC05}, \citeor{Sor09}, \citeor{CGM15} and others, but we do not pursue this direction here.}

Within this framework, our starting point is the following continuous-time \emph{exponential learning} process:
first, each player maintains a vector of \emph{performance scores} that represent his actions' cumulative payoffs;
these scores are then converted into mixed strategies using a logit rule which assigns choice probabilities in proportion to the exponential of each action's score.%
\footnote{The literature refers to versions of this procedure under a variety of names, including weighted/multiplicative majority algorithm (\citeor{LW94}, \citeor{FS99}), exponential weight algorithm (\citeor{Sor09}, \citeor{HSV09}), and Boltzmann Q-learning (\citeor{LC05}, \citeor{THV06}).}
According to a well-known derivation, this logit rule amounts to each player
\revised{maximizing his expected score minus a penalty term}
given by the (negative) entropy of the chosen mixed strategy.
Moreover, under this learning process, mixed strategies evolve according to the replicator dynamics of \citeor{TJ78}, a fundamental model from evolutionary game theory \citeorp{Rus99}.

We extend this framework by considering more general choice maps obtained by replacing the Gibbs entropy with an arbitrary (strongly convex) penalty function that ``regularizes'' the problem of choosing a mixed strategy that maximizes the player's expected score.
These choice maps include those considered in models of stochastic fictitious play (\citeor{FL98}, \citeor{HS02}) and those generated by games with control costs in the refinements literature  (\citeor{vD87}).
In these two cases, the players' choice maps are induced by penalty functions that become infinitely steep at the boundary of the simplex, an assumption which ensures in our context that mixed strategies follow an \ac{ODE} in the interior of the simplex \textendash\ in the case of exponential learning, this is simply the replicator equation of \citeor{TJ78}.
More generally, we find that the evolution of mixed strategies agrees with a class of evolutionary game dynamics studied by \citeor{HS90} and \citeor{Hop99b} (see also \citeor{Har11}), and which we study further in a companion paper \citeorp{MS15}.

Moving beyond these steep cases, our model also allows for penalty functions that are (sub)differentiable over the entire simplex without becoming infinitely steep at the boundary.
The basic example here is the squared Euclidean distance, which induces a choice map based on closest point projection.
In analogue with the logit case, the orbits of this \emph{projected reinforcement learning} process are solutions (in an extended sense) of the projection dynamics of \citeor{Fri91}, another basic model of evolutionary game dynamics.
However, the image of the induced projection-based choice map is the entire simplex (rather than its relative interior), so trajectories of play may now enter and exit the boundary of the simplex in perpetuity.  
Specifically, the \ac{ODE} that describes the evolution of the players' mixed strategies over time only holds for an open dense set of times \textendash\ the union of the open intervals on which the supports of the players' mixed strategies remains constant. 
At all remaining times, the support of the mixed strategy of some player changes, leading to a kink in his mixed strategy trajectory, and the moments at which this occurs cannot always be anticipated by considering the projection dynamics alone.

Our main results extend a variety of basic properties of the replicator dynamics to the class of reinforcement learning dynamics under study, allowing for both steep and nonsteep penalty functions.
First, we show that (iteratively) dominated strategies become extinct and compute their rates of extinction.
Interestingly, while exponential reinforcement learning eliminates dominated strategies at an exponential rate, reinforcement learning with nonsteep penalty functions eliminates such strategies in finite time.
Second, we extend several stability and convergence properties of Nash equilibria under the replicator dynamics.
Namely, we show that
\hspace{-.5ex}\begin{inparaenum}
[\itshape a\upshape)]
\item
limits of (interior) trajectories are \revised{Nash equilibria};
\item
Lyapunov stable states are \revised{Nash equilibria};
and
\item
strict Nash equilibria are asymptotically stable.
\end{inparaenum}
Finally, we show that the basic properties of the time-averaged replicator dynamics (\revised{convergence to equilibrium in zero-sum games that admit an interior equilibrium} and asymptotic agreement with the long-run behavior of the best response dynamics) also extend to the class of learning dynamics studied here.

In the paper most closely related to this one, \citeor{CGM15} (see also \citeor{LC05} and \citeor{THV06}) consider a reinforcement learning process in which players update their mixed strategies based on exponentially discounted payoff estimates.
Focusing exclusively on steep penalty functions, they investigate the convergence of this process in a stochastic, discrete-time environment where players can only observe their realized payoffs.
To do so, the authors also examine the Lyapunov and asymptotic stability properties of (perturbed) Nash equilibria under the resulting mean dynamics in continuous time.
The no-discounting limit of these dynamics coincides with the steep version of the dynamics studied in the current paper, so our stability and convergence results can be seen as an extension of the analysis of \citeor{CGM15} to the nonsteep regime.
Also, \citeor{CGM15} do not examine the elimination of dominated strategies or the long-term behavior of empirical frequencies of play, so our results here provide an indication of other properties that may hold in a discrete-time, stochastic setting.

From the point of view of convex programming, the reinforcement learning dynamics we consider here can be seen as a multi-agent, continuous-time analogue of the well-known \ac{MD} optimization method pioneered by \citeor{NY83} and studied further by \citeor{BecTeb03}, \citeor{ABB04}, \citeor{Nes09} and many others.
This observation also extends to the class of \ac{OMD} algorithms introduced by \citeor{SS11} for \emph{online} convex optimization problems:
focusing on the interplay between discrete- and continuous-time \ac{OMD} schemes, \citeor{KM14} recently showed that a unilateral variant of the reinforcement learning dynamics studied in this paper leads to no regret against any (locally integrable) stream of payoffs.

Our analysis relies heavily on tools from convex analysis and, in particular, the theory of \emph{Bregman functions} \citeorp{Bre67}.
Returning to the case of the replicator dynamics, it is well known that the \emph{\acl{KL} divergence} (an oriented distance measure between probability distributions) is a potent tool for understanding the dynamics' long-run behavior (see \citeor{Wei95} and \citeor{HS98}).
For other steep cases, this part can be played by the \emph{Bregman divergence}, a distance-like function whose role in population dynamics was noted recently by \citeor{Har11}.%
\footnote{Each penalty function \revised{induces} a Bregman divergence;
in the case of the Gibbs entropy, this is simply the \acl{KL} divergence.
For a comprehensive treatment, see \citeor{Kiw97b}.}
In the nonsteep regime however, the dynamics of mixed strategies under reinforcement learning depend intrinsically on the players' score vectors, which determine the supports of the mixed strategies chosen.
To contend with this, we introduce the \emph{Fenchel coupling}, a congruence measure between primal and dual variables
(that is, between mixed strategies and score vectors) that provides a natural tool for proving elimination, convergence, and stability results.

\subsection*{Paper Outline.}

We begin in Section \ref{sec:prelims} with some game-theoretic preliminaries and the definition of penalty functions and choice maps.
The class of reinforcement learning dynamics we examine is presented in Section \ref{sec:dynamics}, along with several examples.
Our analysis proper begins in Section \ref{sec:dominance} where we study the elimination of dominated strategies.
In Section \ref{sec:folk}, we derive some stability and convergence properties of Nash equilibria while, in Section \ref{sec:averages}, we examine the long-term behavior of the players' time-averaged play.
Finally, learning with no regularization is discussed in Section \ref{sec:best}.

\section{Preliminaries}
\label{sec:prelims}

\subsection{Notation.}
\label{sec:notation}

If $\primal$ is a finite-dimensional real space, its dual will be denoted by $\dual$ and we will write $\braket{y}{x}$ for the pairing between $x\in \primal$ and $y\in \dual$.
Also, if $\set = \{s_{\alpha}\}_{\alpha=1}^{n}$ is a finite set, the real vector space \revised{generated} by $\set$ will be denoted by $\R^{\set}$ and its canonical basis by $\{\bvec_{s}\}_{s\in\set}$.
For concision, we use $\alpha$ to refer interchangeably to either $s_{\alpha}$ or $\bvec_{\alpha}$, writing e.g. $x_{\alpha}$ instead of $x_{s_{\alpha}}$;
likewise, we will write $\delta_{\alpha\beta}$ for the Kronecker delta symbols on $\set$.
The set $\simplex(\set)$ of probability measures on $\set$ will be identified with the standard simplex $\simplex = \{x\in \R^{\set}: \sum_{\alpha} x_{\alpha} =1 \text{ and }x_{\alpha}\geq 0\}$ of $\R^{\set}$ and the relative interior of $\simplex$ will be denoted by $\intsimplex \equiv \relint(\simplex)$.
Finally, if $\{\set_{k}\}_{k\in\play}$ is a finite family of finite sets, we will use the shorthand $(\alpha_{k};\alpha_{-k})$ for the tuple $(\dotsc,\alpha_{k-1},\alpha_{k},\alpha_{k+1},\dotsc)$ and we will write $\sum_{\alpha_k}^{k}$ instead of $\sum_{\alpha_k\in\set_{k}}$.

\subsection{Games in normal form.}
\label{sec:games}

A \emph{finite game in normal form} is a tuple $\game \defeq \game(\play,\act,\pay)$ consisting of
\begin{inparaenum}[\itshape a\upshape)]
\item
a finite set of \emph{players} $\play = \{1,\dotsc,N\}$;
\item
a finite set $\act_{k}$ of \emph{actions} (or \emph{pure strategies}) per player $k\in\play$;
and
\item
the players' \emph{payoff functions} $\pay_{k}\from \act\to \R$, where $\act \defeq \prod_{k} \act_{k}$ denotes the game's \emph{action space}, i.e. the set of all \emph{action profiles} $(\alpha_{1},\dotsc,\alpha_{N})$, $\alpha_{k}\in\act_{k}$.
\end{inparaenum}
The set of \emph{mixed strategies} of player $k$ will be denoted by $\strat_{k}\equiv\simplex(\act_{k})$ and the space $\strat\equiv\prod_{k}\strat_{k}$ of \emph{mixed strategy profiles} $x = (x_{1},\dotsc,x_{N})$ will be called the game's \emph{strategy space}.
Unless mentioned otherwise, we will write $\primal_{k} \equiv \R^{\act_{k}}$ for the real space spanned by $\act_{k}$ and $\primal \equiv \prod_{k} \primal_{k}$ for the ambient space of $\strat$.

The expected payoff of player $k$ in the mixed strategy profile $x = (x_{1},\dotsc,x_{N})\in \strat$ is
\begin{equation}
\label{eq:pay}
\pay_{k}(x)
	= \insum_{\alpha_{1}}^{1}\dotsi \insum_{\alpha_{N}}^{N}
	\pay_{k}(\alpha_{1},\dotsc,\alpha_{N}) \; x_{1,\alpha_{1}} \dotsm\, x_{N,\alpha_{N}},
\end{equation}
where $\pay_{k}(\alpha_{1},\dotsc,\alpha_{N})$ denotes the payoff of player $k$ in the profile $(\alpha_{1},\dotsc,\alpha_{N})\in\act$.
Accordingly, the payoff corresponding to $\alpha\in\act_{k}$ in the mixed profile $x\in\strat$ is
\begin{flalign}
\label{eq:payv}
\payv_{k\alpha}(x)
	&\equiv \pay_{k}(\alpha;x_{-k})
	= \insum_{\alpha_{1}'}^{1} \dotsi \insum_{\alpha_{N}'}^{N}
	\pay_{k}(\alpha_{1}',\dotsc,\alpha_{N}')
	\; x_{1,\alpha_{1}'} \dotsm\,\delta_{\alpha,\alpha_{k}'} \dotsm\, x_{N,\alpha_{N}'},
\end{flalign}
and we will write $\payv_{k}(x) = (\payv_{k\alpha}(x))_{\alpha\in\act_{k}}\in \dual_{k}$ for the \emph{payoff \textup(co\textup)vector} of player $k$ at $x\in\strat$.
The prefix ``co'' above is motivated by the natural duality pairing
\begin{equation}
\label{eq:pay-pairing}
\braket{\payv_{k}(x)}{x_{k}}
	= \insum_{\alpha}^{k} x_{k\alpha} \payv_{k\alpha}
	= \pay_{k}(x).
\end{equation}
which shows that $\payv_{k}(x)$ acts on $x_{k}$ as a linear functional.
Duality plays a basic role in our analysis, so mixed strategies $x_{k}\in \primal_{k}$ will be treated as \emph{primal} variables and payoff vectors $\payv_{k}\in \dual_{k}$ as \emph{duals}.%
\footnote{Even though this distinction is rarely made in game theory, it is standard in learning and optimization \textendash\ see e.g. \citeor{Roc70}, \citeor{NY83}, \citeor{SS11} and references therein.}

Finally, a \emph{restriction} of $\game$ is a game $\game'\defeq\game'(\play,\act',u')$ with the same players as $\game$, each with a subset $\act_{k}'\subseteq\act_{k}$ of their original actions and with payoff functions $\pay_{k}'\defeq \pay_{k}|_{\act'}$ suitably restricted to the reduced action space $\act' \eqdef \prod_{k}\act_{k}'$ of $\game'$.

\subsection{Penalty functions and choice maps.}
\label{sec:choice}

In view of \eqref{eq:pay-pairing}, a player's set of optimal mixed strategies given a payoff vector $\payv_{k} \in \dual_{k}$ is 
\begin{equation}
\label{eq:BR}
\brep_{k}(\payv_{k})
	= \argmax_{x_{k}\in\strat_{k}} \braket{\payv_{k}}{x_{k}}.
\end{equation}
A standard way of obtaining a single-valued analogue of the argmax correspondence \eqref{eq:BR} is to introduce a \emph{penalty} term that is (at least) strictly convex in $x_k$.
It is also customary to assume that such penalties are ``infinitely steep'' at the boundary of the simplex (see e.g. \citeor{FL98}), but we obtain a much richer theory by dropping this requirement.
Formally, we have:

\begin{definition}
\label{def:penalty}
Let $\simplex$ be the unit simplex of $\R^{n}$.
We say that $h\from\simplex\to\R$ is a \emph{penalty function} on $\simplex$ if:
\begin{enumerate}

\item
$h$ is continuous on $\simplex$.

\item
$h$ is smooth on the relative interior of every face of $\simplex$ (including $\simplex$ itself).%
\footnote{More precisely, we posit here that $h(\gamma(t))$ is smooth for every smooth curve $\gamma\from(-\eps,\eps)\to\simplex$ that is entirely contained in the relative interior of a given face of $\simplex$.}

\item
$h$ is \emph{strongly convex} on $\simplex$:
there exists some $K>0$ such that
\begin{equation}
\label{eq:str-cvx}
h(tx_{1} + (1-t)x_{2})
	\leq t h(x_{1}) + (1-t) h(x_{2}) - \tfrac{1}{2} K t(1-t) \norm{x_{1}-x_{2}}^{2},
\end{equation}
for all $x_{1}, x_{2}\in\simplex$ and for all $t\in[0,1]$.
\end{enumerate}
If $\norm{dh(x_{n})} \to \infty$ for every interior sequence $x_{n}\to x$, $x_{n}\in\intsimplex$, we will say that $h$ is \emph{steep at} $x$;%
\footnote{In the above, $dh(x)$ for $x\in\intsimplex$ denotes the derivative of the restriction of $h$ to $\intsimplex$, viewed as a map from the tangent space of the simplex to $\R$.
For a detailed treatment, see \citeor[Chap.~3.B.3]{San10}.}
moreover, if this holds for all $x\in\bd(\simplex)$, we will say that $h$ is \emph{steep}.
Finally, $h$ will be called \emph{decomposable with kernel $\theta$} if 
\begin{equation}
\label{eq:decomposable}
h(x)
	= \sum_{\alpha=1}^{n} \theta(x_{\alpha}),
\end{equation}
for some continuous and strongly convex $\theta\from [0,1]\to \R$ that is smooth on $(0,1]$.
\end{definition}

\begin{remark}
\label{rem:fullsteep}
For differentiation purposes, it will often be convenient to assume that the domain of $h$ is the ``thick simplex'' 
$\simplex_\eps = \{x\in \R^{n}: x_{\alpha}\geq0 \text{ and } 1-\eps < \sum_{\alpha} x_{\alpha} < 1 +\eps\}$;
doing so allows us to carry out certain calculations in terms of standard coordinates, but none of our results depend on this device.
Also, ``smooth'' should be interpreted above as ``$C^{\infty}$-smooth'';
our analysis actually requires $C^{2}$ smoothness only if the Hessian of $h$ is involved (and no differentiability otherwise), but we will keep the $C^{\infty}$ assumption for simplicity.
\end{remark}

Given a penalty function $h$ on $\simplex\subseteq \primal \equiv \R^{n}$, the concave maximization problem
\begin{equation}
\label{eq:softmax}
\begin{aligned}
\text{maximize}
	&\quad
	\braket{y}{x} - h(x),
	\\
\text{subject to}
	&\quad
	x\in\simplex,
\end{aligned}
\end{equation}
admits a unique solution for all $y\in \dual$, so the \emph{regularized} correspondence $y\mapsto \argmax_{x} \left\{\braket{y}{x} - h(x)\right\}$ becomes single-valued.
We thus obtain:

\begin{definition}
\label{def:choice}
The \emph{choice map} (or \emph{regularized argmax correspondence}) $\choice\from \dual \to \simplex$ induced by a penalty function $h$ on $\simplex$ is
\begin{equation}
\label{eq:choice}
\choice(y)
	= \argmax\nolimits_{x\in\simplex} \left\{\braket{y}{x} - h(x) \right\},
	\quad
	y\in \dual.
\end{equation}
\end{definition}

Under the steepness and strong convexity requirements of Definition \ref{def:penalty}, the discussion in \citeor[Chapter 26]{Roc70} shows that the induced choice map $\choice$ is smooth and its image is the relative interior $\intsimplex$ of $\simplex$.
At the other end of the spectrum, if $h$ is nowhere steep, the image of $\choice$ is the entire simplex (cf. Remark \ref{rem:imQ} in Appendix \ref{app:Bregman}).%
\footnote{One can also define penalty functions that are steep only at a subset of the boundary of $\simplex$;
the cases considered above are simply the two extremes.}
We illustrate this dichotomy with two representative examples (see also Section \ref{sec:dynamics-examples}):

\begin{example}
\label{ex:choice-logit}
The classic example of a steep penalty function is the (negative) \emph{Gibbs entropy}
\begin{equation}
\label{eq:penalty-Gibbs}
h(x)
	= \sum_{\alpha=1}^{n} x_{\alpha} \log x_{\alpha}.
\end{equation}
As is well known, the induced choice map \eqref{eq:choice} is the so-called \emph{logit map}
\begin{equation}
\label{eq:choice-logit}
\gibbs_\alpha(y)
	= \frac{\exp(y_{\alpha})}{\sum_{\beta=1}^{n} \exp(y_{\beta})} .
\end{equation} 
Since $h$ is steep, $\im \gibbs = \intsimplex$.
\end{example}

\begin{example}
\label{ex:choice-Eucl}
The standard example of a non-steep penalty function is the \emph{quadratic penalty}
\begin{equation}
\label{eq:penalty-Eucl}
h(x)
	= \frac{1}{2} \sum_{\alpha=1}^{n} x_{\alpha}^{2}.
\end{equation}
The induced choice map \eqref{eq:choice} is the (Euclidean) \emph{projection map}
\begin{equation}
\label{eq:choice-Eucl}
\txs
\eucl(y)
	= \argmax_{x\in\simplex} \left\{ \braket{y}{x} - \tfrac{1}{2}\norm{x}_{2}^{2} \right\}
	= \argmin_{x\in\simplex}
	\norm{y - x}_{2}^{2}
	= \proj_{\simplex} y,
\end{equation}
where $\proj_{\simplex}$ denotes the closest point projection to $\simplex$ with respect to the standard Euclidean norm $\norm{\argdot}_{2}$ on $\R^{n}$.  Obviously, $\im \eucl = \simplex$.
\end{example}


\begin{remark}
Up to mild technical differences, penalty functions are also known as \emph{regularizers} in online learning and \emph{Bregman functions} (or \emph{prox-functions}) in convex analysis;
for a comprehensive treatment, see \citeor{Bre67}, \citeor{NY83}, \citeor{SS11} and references therein.
The term ``decomposable'' is borrowed from \citeor{ABB04}.

In game theory, \citeor{FL98} use the term \emph{smooth best response function} to refer to the composition $\choice_{k} \circ \payv_{k}$ of a choice map $\choice_{k}$ generated by a steep penalty function and a game's payoff function $\payv_{k}$.  
\citeor{HS02} use the term ``perturbed best response function'' for such composite functions, while \citeor{MP95} use the term ``quantal response function'' to refer to $\choice_{k}$ directly (the notation $\choice$ is in reference to this last fact).
\end{remark}

\begin{remark}
\label{rem:StrongConvexity}
Several models of smooth fictitious play (\citeor{FL98}, \citeor{HS02}) use a steep penalty function with positive-definite Hessian.
Concerning this last condition, the strong convexity of $h$ imposes a positive lower bound on the smallest eigenvalue of $\hess(h)$, a property which in turn ensures that the associated choice map $\choice$ is Lipschitz continuous (Proposition \ref{prop:choice-dh}).
Later, we also take advantage of the fact that strong convexity provides a lower bound for the so-called \emph{Bregman divergence} between points in $\simplex$ (Proposition \ref{prop:Bregman}).
Even though some of our results can be extended to penalty functions that are not strongly convex (for instance, decomposable penalty functions with a strictly convex kernel), the above consequences of strong convexity simplify our presentation considerably so we will not venture beyond the strongly convex case.
\end{remark}

\section{A class of reinforcement learning dynamics}
\label{sec:dynamics}

\subsection{Definition and basic examples.}
\label{sec:RL}

The basic reinforcement learning scheme that we consider is that players
keep track of the cumulative payoffs of their actions and then use a choice map to transform these aggregate score vectors into mixed strategies and keep playing.
More precisely, given a finite game $\game \equiv \game(\play,\act,\pay)$, we will focus on the continuous-time process
\begin{equation}
\label{eq:RL}
\tag{RL}
\begin{aligned}
\score_{k}(t)
	&= \score_{k}(0) + \int_{0}^{t} \payv_{k}(x(s)) \dd s,
	\\
x_{k}(t)
	&= \choice_{k}(\score_{k}(t)),
\end{aligned}
\end{equation}
or, in differential form:
\begin{equation}
\label{eq:RL-diff}
\dot \score_{k}
	= \payv_{k}(\choice(\score)),
\end{equation}
where $\choice \equiv (\choice_{1},\dotsc,\choice_{N})\from \dual \equiv \prod_{k} \dual_{k}\to\strat$ and $\score = (\score_{1},\dotsc,\score_{N})\in \dual$ denote the players' choice and score profiles respectively.
In the above, the (primal) strategy variable $x_{k}(t) \in \strat_{k}$ describes the mixed strategy of player $k$ at time $t$ while the (dual) score vector $\score_{k}(t)$ aggregates the payoffs of the pure strategies $\alpha\in\act_{k}$ of player $k$.
Accordingly, the basic interpretation of \eqref{eq:RL} is that each player observes the realized expected payoffs of his strategies over a short interval of time and then uses these payoffs to update his score vector.

\newcommand{\dtime}{n}

More precisely, in the stochastic approximation language of \citeor{Ben99}, \eqref{eq:RL} is simply the mean field of the discrete-time stochastic process
\begin{equation}
\label{eq:RL-discrete}
\begin{aligned}
\Score_{k\alpha}(\dtime+1)
	&= \Score_{k\alpha}(\dtime)
	+ \hat\payv_{k\alpha}(\dtime),
	\\
X_{k\alpha}(\dtime+1)
	&= \choice_{k\alpha}(\Score_{k}(\dtime+1)),
\end{aligned}
\end{equation}
where $X_{k\alpha}(\dtime)$ is the probability of playing $\alpha\in\act_{k}$ at the $\dtime$-th instance of play, $\dtime=0,1,\dotsc$, while $\hat\payv_{k\alpha}(\dtime)$ is an unbiased estimator of $\payv_{k\alpha}(X(\dtime))$.
If player $k$ can observe the action profile $\alpha_{-k}(\dtime)$ played by his opponents (or can otherwise calculate his strategies' payoffs), such an estimate is provided by $\hat\payv_{k\alpha}(\dtime) = \pay_{k}(\alpha;\alpha_{-k}(\dtime))$.
Instead, if player $k$ can only observe the payoff $\hat\pay_{k}(\dtime) = \pay_{k}(\alpha_{k}(\dtime);\alpha_{-k}(\dtime))$ of his chosen action $\alpha_{k}(\dtime)$, a standard choice for $\hat\payv_{k}(\dtime)$ is 
\begin{equation}
\label{eq:scoredivx}
\hat\payv_{k\alpha}(\dtime)
	= 
	\hat\pay_{k}(\dtime) / X_{k\alpha}(\dtime)
	\quad
	\text{if $\alpha_{k}(\dtime) = \alpha$},
\end{equation}
where division by $X_{k\alpha}(\dtime)$ compensates for the infrequency with which the score of strategy $\alpha$ is updated.
Estimator \eqref{eq:scoredivx} is sound if the penalty function of player $k$ is steep (see \citeor{LC05} and \citeor{CGM15});
otherwise, $X_{k\alpha}(\dtime)$ may become zero, in which case the links between \eqref{eq:scoredivx} and our continuous-time model are less clear.

We begin with two representative examples of the reinforcement learning dynamics \eqref{eq:RL}:

\begin{example}
\label{ex:ExpRIP}
If $\choice$ is the logit map \eqref{eq:choice-logit} of Example \ref{ex:choice-logit}, players select actions with probability proportional to the exponential of their aggregate payoffs.
In this case, \eqref{eq:RL} boils down to the \emph{exponential} (or \emph{logit}) \emph{reinforcement learning process}
\begin{equation}
\label{eq:XL}
\tag{XL}
\begin{aligned}
\dot\score_{k\alpha}
	&= \payv_{k\alpha}(x),\\
x_{k\alpha}
	&= \frac{\exp(\score_{k\alpha})}{\insum_{\beta}^{k} \exp(\score_{k\beta})}.
\end{aligned}
\end{equation}
In a single-agent, online learning context, the discrete-time version of \eqref{eq:XL} first appeared in the work of \citeor{Vov90} and \citeor{LW94} \textendash\ see also \citeor{Rus99}, \citeor{Sor09}, and \citeor{KM14} for a continuous-time analysis.
In a game-theoretic setting, this process has been studied by (among others) \citeor{FS99}, \citeor{HSV09} and \citeor{MM10}, while \citeor{LC05}, \citeor{THV06} and, more recently, \citeor{CGM15} considered a discounted variant that we describe in Section \ref{sec:related} below.

Differentiating $x_{k\alpha}$ in \eqref{eq:XL} with respect to time and substituting yields 
\begin{equation}
\label{eq:RD-from-XL}
\dot x_{k\alpha}
	= \frac
	{\dot \score_{k\alpha} e^{\score_{k\alpha}} \insum_{\beta}^{k} e^{\score_{k\beta}}
	- e^{\score_{k\alpha}} \insum_{\beta}^{k} \dot \score_{k\beta} e^{\score_{k\beta}}}
	{\left(\insum_{\beta}^{k} e^{\score_{k\beta}} \right)^{2}}
	= x_{k\alpha} \left[ \dot\score_{k\alpha} - \insum_{\beta}^{k} x_{k\beta} \dot\score_{k\beta} \right],
\end{equation}
so, with $\dot \score_{k\alpha} = \payv_{k\alpha}(x)$, we readily obtain:
\begin{equation}
\label{eq:RD}
\tag{RD}
\dot x_{k\alpha}
	= x_{k\alpha} \left[\payv_{k\alpha}(x) - \insum_{\beta}^{k} x_{k\beta} \payv_{k\beta}(x) \right].
\end{equation}
This equation describes the replicator dynamics of \citeor{TJ78},
a fundamental model of evolutionary game theory whose long-term rationality properties are quite well understood.
This basic relation between exponential reinforcement learning and the replicator dynamics was noted in a single-agent environment by \citeor{Rus99} and was explored further in a game-theoretic context by \citeor{HSV09} and \citeor{MM10}.
\end{example}

\begin{figure}[t]
\centering
\includegraphics[width=.95\textwidth]{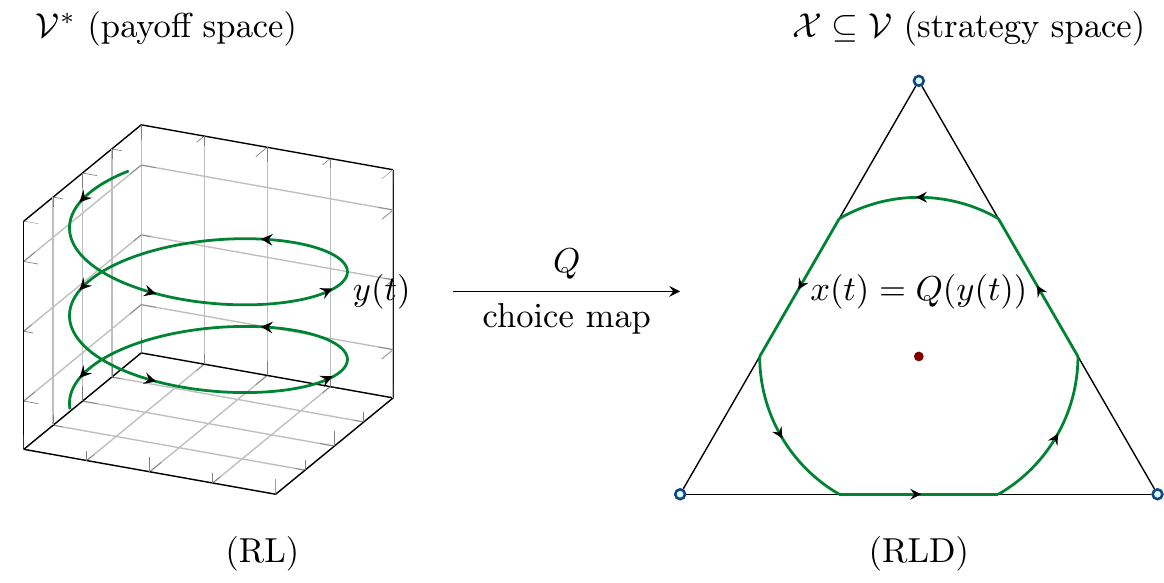}
\smallskip
\caption{\small
Diagrammatic illustration of the relation between the reinforcement learning scheme \eqref{eq:RL} and the induced dynamics \eqref{eq:RLD}.
The dynamics \eqref{eq:RL} evolve in the dual space $\dual$ and solutions exist for all time;
on the other hand, the dynamics \eqref{eq:RLD} evolve in the primal space $\strat$ and solutions may fail to exist for all time.
The solution orbits $x(t) = \choice(y(t))$ of \eqref{eq:RL} only fail to satisfy \eqref{eq:RLD} when they pass from one face of $\strat$ to another.}
\label{fig:choice}
\end{figure}

\begin{example}
If $\choice$ is the projection map \eqref{eq:choice-Eucl} of Example \ref{ex:choice-Eucl}, \eqref{eq:RL} leads to the \emph{projected reinforcement learning} process
\begin{equation}
\label{eq:PL}
\tag{PL}
\begin{aligned}
\dot\score_{k}
	&= \payv_{k}(x),
	\\
x
	&= \proj_{\strat} \score.
\end{aligned}
\end{equation}
Of course, since $\proj_{\strat} y$ is not smooth in $y$, we can no longer use the same approach as in \eqref{eq:RD-from-XL} to derive the dynamics of the players' mixed strategies $x_{k}$.
Instead, recall (or solve the defining convex program to show) that the closest point projection on $\strat_{k} = \simplex(\act_{k})$ takes the simple form
\begin{equation}
\label{eq:proj-simplex}
\mathopen\big( \proj_{\strat_{k}} y_{k} \big)_{\alpha}
	=\max\{y_{k\alpha} + \mu_{k},0\},
\end{equation}
where $\mu_{k}\in\R$ is such that $\sum_{\alpha}^{k} \max\{y_{k\alpha} + \mu_{k},0\} = 1$.
Therefore, if $\interval_{k}$ is an open time interval over which $x_{k}(t) = \proj_{\strat_{k}} \score_{k}(t)$ has constant support $\act_{k}' \subseteq\act_{k}$,
a simple differentiation yields
\begin{equation}
\label{eq:xdot-Eucl}
\dot x_{k\alpha}
	= \dot \score_{k\alpha} + \dot \mu_{k}
	= \payv_{k\alpha} + \dot\mu_{k}
	\quad
	\text{for all $\alpha\in\act_{k}'$.}
\end{equation}
Since $\sum_{\alpha\in\act_{k}'} \dot x_{k\alpha} = 0$, summing over $\alpha\in\act_{k}'$ gives
\begin{equation}
0
	= \insum_{\alpha\in\act_{k}'} \payv_{k\alpha} + \dot\mu_{k} \abs{\act_{k}'},
\end{equation}
so, by substituting into \eqref{eq:xdot-Eucl} and rearranging, we obtain the \emph{projection dynamics}:
\begin{equation}
\label{eq:PD}
\tag{PD}
\dot x_{k\alpha}
	= \begin{cases}
	\dis\payv_{k\alpha}(x) - \abs{\supp(x_{k})}^{-1} \insum_{\beta\in\supp(x_{k})} \payv_{k\beta}(x)
		&\text{if $\alpha\in\supp(x_{k})$,}
		\\
	0
		&\text{if $\alpha\notin \supp(x_{k})$.}
	\end{cases}
\end{equation}

The dynamics \eqref{eq:PD} were introduced in game theory by \citeor{Fri91} as a geometric model of the evolution of play in population games.%
\footnote{\citeor{NZ97} (see also \citeor{LS08} and \citeor{SDL08}) introduce related projection-based dynamics for population games.
The relations among the various projection dynamics are explored in a companion paper \citeorp{MS15}.}
The previous discussion shows that the projected orbits $x(t) = \proj_{\strat} y(t)$ of the learning scheme \eqref{eq:PL} satisfy the projection dynamics \eqref{eq:PD} on every open interval over which the support of $x(t)$ is fixed;
furthermore, as we argue below, the union of these intervals is dense in $[0, \infty)$.
In this way, orbits of \eqref{eq:PL} that begin in the relative interior $\intstrat$ of the game's strategy space may attain a boundary face in finite time, then move to another boundary face or re-enter $\intstrat$ (again in finite time), and so on (cf. Fig.~\ref{fig:choice}).
Thus, although $x(t)$ may fail to be differentiable when it moves from (the relative interior of) one face of $\strat$ to another, it satisfies \eqref{eq:PD} for all times in between.
\end{example}

These two examples illustrate a fundamental dichotomy between reinforcement learning processes induced by steep and nonsteep penalty functions.
In the steep case, the dynamics of the strategy variable are well-posed and admit unique solutions that stay in $\intstrat$ for all time.
On the other hand, in the nonsteep regime, the dynamics of the strategy variable only admit solutions in an extended sense, and they may enter or exit different faces of $\strat$ in perpetuity.  

\begin{remark}
In several treatments of stochastic fictitious play (\citeor{FL98}, \citeor{HS02}), it is common to replace $h(x)$ with $\temp h(x)$ for some positive parameter $\temp > 0$ which is often called the model's \emph{noise level}.%
\footnote{The term ``noise level'' reflects the fact that $\temp$ essentially controls the magnitude of the perturbation to the player's expected payoff in the regularized maximization problem \eqref{eq:softmax}.}
For instance, if $h$ is the Gibbs entropy \eqref{eq:penalty-Gibbs}, this leads to the choice map
\begin{equation}
\label{eq:choice-logit2}
\gibbs_{\alpha}^{\tempinv}(y)
	= \frac{\exp(\tempinv y_{\alpha})}{\sum_{\beta=1}^{n} \exp(\tempinv y_{\beta})} .
\end{equation} 
with $\tempinv = \temp^{-1}$.
If $y$ is fixed, choices are nearly uniform for small $\tempinv$;
on the other hand, for large $\tempinv$, almost all probability is placed on the pure strategies with the highest score.

In the present context, replacing $h_k(x_{k})$ with $\temp_{k} h_{k}(x_{k})$ and writing $\tempinv_{k}= \temp_{k}^{-1}$ yields the following variant of \eqref{eq:RL}:
 \begin{equation}
\label{eq:RL-temp}
\begin{aligned}
\dot \score_{k}
	&=  \payv_{k}(x),
	\\
x_{k}
	&= \choice_{k}(\tempinv_{k}\score_{k}).
\end{aligned}
\end{equation}
Since a rescaled penalty function is still a penalty function, \eqref{eq:RL-temp} can be viewed as an instance of \eqref{eq:RL}; therefore, our results for the latter also apply to the former.
Furthermore, because the score variables $\score_k(t)$ scale with $t$, introducing $\tempinv$ has less drastic consequences under \eqref{eq:RL} than under stochastic fictitious play:
for instance, the stationary points of \eqref{eq:RL} in $\strat$ remain unaffected by this choice \textendash\ see Theorem \ref{thm:folk} below.

One can also consider a variant of \eqref{eq:RL} under which different players adjust their score variables at different rates:
\begin{equation}
\label{eq:RL-rate}
\tag{\ref*{eq:RL}$_{\tempinv}$}
\begin{aligned}
\dot \score_{k}
	&= \tempinv_{k} \payv_{k}(x),
	\\
x_{k}
	&= \choice_{k}(\score_{k}),
\end{aligned}
\end{equation}
Evidently, this process is equivalent to \eqref{eq:RL-temp}, but with initial conditions $y_k(0)$ scaled by $1/\tempinv_{k}$.
As we show in Propositions \ref{prop:dom-rate} and \ref{prop:strict-rate}, the choice of $\tempinv$ affects the speed at which \eqref{eq:RL-temp} evolves because it determines each player's characteristic time scale.
\end{remark}

\subsection{Related models.}
\label{sec:related}

Before proceeding with our analysis of \eqref{eq:RL}, we mention a number of related models appearing in the literature.

First, as an alternative to aggregating payoffs in \eqref{eq:RL}, one can consider the exponentially discounted model
\begin{equation}
\label{eq:discount}
\score_{k}(t)
	= \score_{k}(0) \lambda^{t} + \int_{0}^{t} \lambda^{t-s} \payv_{k}(x(s)) \dd s,
\end{equation}
where the discount rate $\lambda\in(0,1)$ measures the relative weight of past observations.
This variant was examined by \citeor{LC05}, \citeor{THV06}, and \citeor{CGM15} for choice maps generated by steep penalty functions.
Obviously, when $h$ is steep, \eqref{eq:RL} can be seen as a limiting case of \eqref{eq:discount} for $\lambda\to1^{-}$.
In contrast to \eqref{eq:RL} however, discounting implies that the score variable  $\score_{k}(t)$ remains bounded, thus preventing the agents' mixed strategies from approaching the boundary of $\strat$.
For instance, under the logit rule of Example \ref{ex:ExpRIP}, \citeor{CGM15} showed that discounting introduces a penalty term which repels orbits from the boundary $\bd(\strat)$ of $\strat$ under the replicator dynamics.

In a single-agent environment, payoffs are determined at each instance by nature so the reinforcement learning process \eqref{eq:RL} becomes
\begin{equation}
\label{eq:RL-unilateral}
\score_{\alpha}(t)
	= \score_{\alpha}(0)
	+ \int_{0}^{t} \payv_{\alpha}(s) \dd s,
	\qquad
x(t)
	= \choice(\score(t)).
\end{equation}
In this context, \eqref{eq:RL} can be seen as a continuous-time analogue of the family of online learning algorithms known as \acf{OMD} \textendash\ for a comprehensive account, see \citeor{Bub11} and \citeor{SS11}.
The resulting interplay between discrete and continuous time has been analyzed by \citeor{Sor09} and \citeor{KM14} who showed that \eqref{eq:RL-unilateral} leads to no regret against any locally integrable payoff stream $\payv(t)$ in $\dual$.
In view of the above, \eqref{eq:RL} extends the discounted dynamics of \citeor{LC05} to the nonsteep regime and the online learning dynamics of \citeor{KM14} to a multi-agent, game-theoretic setting.

There are several other reinforcement learning schemes that are distinct from \eqref{eq:RL} but which still lead to the replicator equation \eqref{eq:RD}.
A leading model of this kind is presented in the seminal paper of \citeor{ER98}:
after player $k$ chooses pure strategy $\alpha\in\act_{k}$, he increments its score by the payoff he receives (assumed positive) and then updates his choice probabilities proportionally to each action's score.
The continuous-time, deterministic version of this model is
\begin{equation}
\label{eq:ER}
\tag{ER}
\begin{aligned}
\dot\score_{k\alpha}
	&= x_{k\alpha}\payv_{k\alpha}(x),\\
x_{k\alpha}
	&= \frac{\score_{k\alpha}}{\insum_{\beta}^{k} \score_{k\beta}},
\end{aligned}
\end{equation}
where the $x_{k\alpha}$ term in the first equation reflects the fact that the score variable $\score_{k\alpha}$ is only updated when $\alpha$ is played.
A simple calculation then shows that the evolution of mixed strategies is governed by the replicator dynamics \eqref{eq:RD} up to a player-specific multiplicative factor.
Versions of this model have been studied by \citeor{Pos97}, \citeor{Rus99}, \citeor{Hop02}, \citeor{Beg05}, and \citeor{HP05}; 
\citeor{Rus99} also considers hybrids between \eqref{eq:XL} and \eqref{eq:ER} in nonstrategic settings.

\citeor{BS97} also consider a variant of the learning model of \citeor{Cro73} where there is no separate score variable.
Instead, if player $k$ chooses pure strategy $\alpha$, he increases the probability with which he plays $\alpha$ and decreases the probability of every other action $\beta\neq\alpha$ proportionally to the payoff $\payv_{k\alpha}(x)\in(0,1)$ that the player obtained.  
The continuous-time, deterministic version of this model is
\begin{equation}
\dot x_{k\alpha}
	= x_{k\alpha}(1- x_{k\alpha}) \payv_{k\alpha}(x)
	+ \insum_{\beta \neq \alpha}x_{k\beta} (0- x_{k\alpha}) \payv_{k\beta}(x),
\end{equation}
which again yields the replicator dynamics \eqref{eq:RD} after a trivial rearrangement.

\subsection{Basic results.}
\label{sec:basic}

We start our analysis by showing that the dynamics \eqref{eq:RL} are well-posed even if the players' penalty functions are not steep:

\begin{proposition}
\label{prop:wp}
The reinforcement learning process \eqref{eq:RL} admits a unique global solution for every initial score profile $y(0)\in \dual$.
\end{proposition}

\begin{Proof}
By Proposition \ref{prop:choice-dh}, strong convexity implies that $\choice\from \dual \to \strat$ is Lipschitz.
Since $\payv_{k}$ is bounded, existence and uniqueness of global solutions follows from standard arguments \textendash\ e.g.~\citeor[Chapter~V]{Rob95}.
\end{Proof}

We turn now to the dynamics induced by \eqref{eq:RL} on the game's strategy space $\strat$.
To that end, if $y(t)$ is a solution orbit of \eqref{eq:RL}, we call $x(t) = \choice(y(t))$ the \emph{trajectory of play induced by $y(t)$} \textendash\ or, more simply, an \emph{orbit of \eqref{eq:RL} in $\strat$}.
Mirroring the derivation of the projection dynamics \eqref{eq:PD} above, our next result provides a dynamical system on $\strat$ that is satisfied by smooth segments of orbits of \eqref{eq:RL} in $\strat$.
To state it, let
\begin{equation}
g_{k}^{\alpha\beta}(x_{k})
	= \hess\big(h_{k}\vert_{\strat_{k}'}(x_{k})\big)_{\alpha\beta}^{-1},
	\qquad
	\alpha,\beta\in\supp(x_{k}),
\end{equation}
denote the inverse Hessian matrix of the restriction $h_{k}\vert_{\strat_{k}'}$ of $h_{k}$ to the face $\strat_{k}' = \simplex(\supp(x_{k}))$ of $\strat_{k}$ that is spanned by $\supp(x_{k})$.%
\footnote{Strong convexity ensures that \revised{$\hess(h_{k}\vert_{\strat_{k}'})$} is positive-definite \textendash\ and, hence, invertible (cf. Remark \ref{rem:StrongConvexity}).}
Furthermore, let
\begin{equation}
g_{k}^{\alpha}(x_{k})
	= \insum_{\beta\in\supp(x_{k})} g_{k}^{\alpha\beta}(x_{k})
	\text{\quad and \quad}
G_{k}(x_{k})
	 = \insum_{\alpha\in\supp(x_{k})} g_{k}^{\alpha}(x_{k})
\end{equation}
denote the row sums and the grand sum of $g_{k}^{\alpha\beta}(x_{k})$ respectively.
We then have:

\begin{proposition}
\label{prop:RLD}
Let $x(t) = \choice(y(t))$ be an orbit of \eqref{eq:RL} in $\strat$, and let $\interval$ be an open interval over which the support of $x(t)$ remains constant.
Then, for all $t\in\interval$, $x(t)$ satisfies:
\begin{equation}
\label{eq:RLD}
\tag{\ref*{eq:RL}D}
\dot x_{k\alpha}
	= \sum_{\beta\in\supp(x_{k})}
	\left[ g_{k}^{\alpha\beta}(x_{k}) - \frac{1}{G_{k}(x_{k})} g_{k}^{\alpha}(x_{k}) g_{k}^{\beta}(x_{k}) \right] \payv_{k\beta}(x),
	\quad
	\alpha\in\supp(x_{k}).
\end{equation}
\end{proposition}

\begin{corollary}
\label{cor:RLD-steep}
Every orbit $x(t) = \choice(y(t))$ of \eqref{eq:RL} in $\strat$ is Lipschitz continuous and satisfies \eqref{eq:RLD} on an open dense subset of $[0,\infty)$. 
Furthermore, if the players' penalty functions are steep, \eqref{eq:RLD} is well-posed and $x(t)$ is an ordinary solution thereof.
\end{corollary}

\begin{Proof}[Proof of Proposition \ref{prop:RLD}]
Let $\act_{k}'$ denote the (constant) support of $x_{k}(t)$ for $t\in\interval$.
Then, the first-order \ac{KKT} conditions for the softmax problem \eqref{eq:softmax} of player $k$ readily yield
\begin{equation}
\label{eq:KKT-interior}
\score_{k\alpha} - h_{k\alpha}
	= \mu_{k}
	\quad
	\text{for all $\alpha\in\act_{k}'$,}
\end{equation}
where $\mu_{k}$ is the Lagrange multiplier associated to the constraint $\sum_{\alpha\in\act_{k}'} x_{k\alpha} = 1$ and we have set $h_{k\alpha} = \frac{\pd}{\pd x_{k\alpha}} h_{k}\vert_{\strat_{k}'}$.
Differentiating \eqref{eq:KKT-interior} then yields
\begin{equation}
\label{eq:der1}
\dot \score_{k\alpha}
	- \sum_{\beta\in\act_{k}'} \hess(h_{k}\vert_{\strat_{k}'}(x))_{\alpha\beta}\, \dot x_{k\beta}
	= \dot \mu_{k},
\end{equation}
and hence
\begin{equation}
\label{eq:der2}
\dot x_{k\alpha}
	= \sum_{\beta\in\act_{k}'} g_{k}^{\alpha\beta}(x_{k}) \left( \payv_{k\beta}(x) - \dot \mu_{k} \right)
	= \sum_{\beta\in\act_{k}'} g_{k}^{\alpha\beta}(x_{k}) \payv_{k\beta}(x) - g_{k}^{\alpha}(x_{k}) \dot \mu_{k},
\end{equation}
where we have used the fact that $\dot \score_{k} = \payv_{k}$.
However, since $x_{k}(t) \in \strat_{k}'$ for all $t\in\interval$ by assumption, we must also have $\sum_{\alpha\in\act_{k}'} \dot x_{k\alpha} = 0$;
accordingly, \eqref{eq:der2} gives
\begin{equation}
\label{eq:der3}
0
	= \sum_{\alpha,\beta\in\act_{k}'} g_{k}^{\alpha\beta}(x_{k}) \payv_{k\beta}(x)
	- \dot \mu_{k} \sum_{\alpha\in\act_{k}'} g_{k}^{\alpha}(x_{k})
	= \sum_{\beta\in\act_{k}'} g_{k}^{\beta}(x_{k}) \payv_{k\beta}(x)
	- G_{k}(x_{k}) \dot \mu_{k},
\end{equation}
so \eqref{eq:RLD} is obtained by solving \eqref{eq:der3} for $\dot\mu_{k}$ and substituting in \eqref{eq:der2}.
\end{Proof}

\begin{Proof}[Proof of Corollary \ref{cor:RLD-steep}]
Lipschitz continuity follows from Proposition \ref{prop:wp} and the Lipschitz continuity of $\choice$ (Proposition \ref{prop:choice-dh}).
To establish the next claim, we must show that the union of all open intervals over which $x(t)$ has constant support is dense in $[0,\infty)$.
To do so, fix some $\alpha\in\act_{k}$, $k\in\play$, and let $A = \{t: x_{k\alpha}(t) > 0\}$ so that $A^{c} = x_{k\alpha}^{-1}(0)$.
Then, if $B=\inter(A^{c})$, it suffices to show that $A\cup B$ is dense in $[0, \infty)$.
Indeed, if $t \notin \cl(A \cup B)$, we must have $t \notin A$ and hence $x_{k\alpha}(t) = 0$. 
Furthermore, since $t \notin \cl(A)$, there exists a neighborhood $U$ of $t$ that is disjoint from $A$, i.e. $x_{k\alpha} = 0$ on $U$.
Since $U$ is open, we get $U \subseteq \inter(A^{c}) = B$, contradicting that $t \notin B$.

Finally, to prove the second part of our claim, simply note that the image $\im \choice_{k}$ of $\choice_{k}$ coincides with the relative interior $\intstrat_{k}$ of $\strat_{k}$ if and only if $h_{k}$ is steep (cf. Proposition \ref{prop:choice-dh} in Appendix \ref{app:Bregman}).
\end{Proof}

\begin{remark}
If the players' penalty functions can be decomposed as $h_{k}(x_{k}) = \insum_{\beta}^{k} \theta_{k}(x_{k\beta})$ (cf. Definition \ref{def:penalty}), the inverse Hessian matrix of $h_{k}$ may be written as
\begin{equation}
\label{eq:ginv-theta}
g_{k}^{\alpha\beta}(x_{k})
	= \frac{\delta_{\alpha\beta}}{\theta_{k}''(x_{k\alpha})},
	\qquad
	\alpha,\beta\in\supp(x_{k}).
\end{equation}
In this case, \eqref{eq:RLD} may be written more explicitly as
\begin{equation}
\label{eq:RLD-theta}
\tag{\ref*{eq:RLD}$_{\theta}$}
\dot x_{k\alpha}
	= \frac{1}{\theta_{k}''(x_{k\alpha})}
	\left[
	\payv_{k\alpha}(x) - \Theta_{k}''(x_{k}) \insum_{\beta\in\supp(x_{k})} \payv_{k\beta}(x)\big/\theta_{k}''(x_{k\beta})
	\right],
\end{equation}
where $\Theta_{k}''$ stands for the harmonic aggregate%
\footnote{We should stress here that $\Theta_{k}''$ is not a second derivative;
we only use this notation for visual consistency.}
\begin{equation}
\label{eq:Theta}
\Theta_{k}''(x_{k})
	= \left[\insum_{\beta\in\supp(x_{k})} 1/\theta_{k}''(x_{k\beta})\right]^{-1}.
\end{equation}

To the best of our knowledge, the dynamics \eqref{eq:RLD-theta} first appeared in a comparable form in the work of \citeor{Har11} under the name ``escort replicator dynamics''.%
\footnote{See also \citeor{CGM15} for a variant of \eqref{eq:RLD-theta} induced by the exponentially discounted model \eqref{eq:discount}.}
From the perspective of convex programming, the dynamics \eqref{eq:RLD-theta} for steep $\theta$ (and, more generally, \eqref{eq:RLD} for steep $h$) can be seen as a game-theoretic analogue of the Hessian Riemannian gradient flow framework of \citeor{BT03} and \citeor{ABB04}.
As such, \eqref{eq:RLD} exhibits a deep Riemannian-geometric character which links it to class of dynamics introduced by \citeor{HS90} and studied further by \citeor{Hop99a}.
These geometric aspects of \eqref{eq:RLD} are explored in detail in a companion paper \citeorp{MS15}.
\end{remark}

\begin{remark}
One subtle point in Corollary \ref{cor:RLD-steep} is that $x(t) = \choice(y(t))$ need not be a solution of \eqref{eq:RLD} in the sense of Carathéodory.
The reason for this is that Carathéodory solutions are required to satisfy the dynamical system at hand over a set of full measure;
by contrast, Corollary \ref{cor:RLD-steep} shows that $x(t)$ satisfies \eqref{eq:RLD} over an \emph{open dense} subset of times.
Hence, in principle, $x(t)$ may fail to satisfy \eqref{eq:RLD} over a closed, nowhere dense set with \emph{positive} measure \textendash\ such as a fat Cantor set.
We believe that intricate topological pathologies of this sort do occur under \eqref{eq:RLD}, but we have not been able to prove it either.
\end{remark}


\subsection{Further examples.}
\label{sec:dynamics-examples}

In Section \ref{sec:RL}, we introduced the exponential and projected reinforcement learning models \eqref{eq:XL} and \eqref{eq:PL}, generated respectively by the (steep) entropic penalty \eqref{eq:penalty-Gibbs} and the (nonsteep) quadratic penalty \eqref{eq:penalty-Eucl}.
We close this section with some further examples of penalty functions and their induced mixed strategy dynamics.

\begin{example}[The Tsallis entropy and the $q$-replicator dynamics]
\label{ex:qRD}

A well known generalization of the Gibbs (negative) entropy due to \citeor{Tsa88} is:%
\footnote{In information theory, the Tsallis entropy is often referred to as the Havrda\textendash Charvát entropy.
Also, note that \eqref{eq:penalty-Tsallis} uses the normalization $[q(1-q)]^{-1}$ rather than the more common $(1-q)^{-1}$;
this is done to simplify notation later on.}
\begin{equation}
\label{eq:penalty-Tsallis}
h(x)
	= [q(1-q)]^{-1}\sum_{\alpha=1}^{n} (x_{\alpha} - x_{\alpha}^{q}),
	\qquad
	q>0,
\end{equation}
with the continuity convention $(z - z^{q})/(1-q) = z \log z$ for $q=1$ (corresponding to the Gibbs penalty of Example \ref{ex:choice-logit}).
This penalty function is decomposable (in the sense of Definition \ref{def:penalty}) with kernel $\theta(x) = [q(1-q)]^{-1} (x - x^{q})$. A simple differentiation then gives $\theta''(x) = x^{q-2}$ and $\Theta''(x) \equiv \big[\insum_{\beta} 1/\theta''(x_{\beta}) \big]^{-1} = \big[ \insum_{\beta} x_{\beta}^{2-q} \big]^{-1}$ so, substituting in \eqref{eq:RLD-theta}, we obtain the \emph{$q$-replicator dynamics}
\begin{equation}
\label{eq:qRD}
\tag{RD$_{q}$}
\dot x_{k\alpha}
	= x_{k\alpha}^{2-q} \left[\payv_{k\alpha} - r_{k,q}^{q-2}\insum_{\beta}^{k} x_{k\beta}^{2-q} \payv_{k\beta}\right],
\end{equation}
where we set $r_{k,q}^{2-q} = \insum_{\beta}^{k} x_{k\beta}^{2-q}$ (and we are using the convention $0^{0} = 0$ for $q=2$).

In the context of convex programming, \eqref{eq:qRD} was derived as an example of a \acl{HR} gradient flow in \citeor{ABB04}.
More recently, these dynamics also appeared in \citeor{Har11} under the name ``$q$-deformed replicator dynamics''.
Obviously, for $q=1$, \eqref{eq:qRD} is simply the replicator equation \eqref{eq:RD}, reflecting the fact that the Tsallis penalty \eqref{eq:penalty-Tsallis} converges to the Gibbs penalty \eqref{eq:penalty-Gibbs} as $q\to1$.
Furthermore, for $q=2$, the Tsallis penalty \eqref{eq:penalty-Tsallis} is equal to the quadratic penalty \eqref{eq:penalty-Eucl} up to an affine term;
consequently, since \eqref{eq:RLD-theta} does not involve the first derivatives of $h$, the dynamics \eqref{eq:qRD} for $q=2$ are the same as the projection dynamics \eqref{eq:PD}.

Of course, \eqref{eq:penalty-Tsallis} is steep if and only if $0<q\leq1$, so \eqref{eq:qRD} may fail to be well-posed for $q>1$.
In particular, as in the case of the projection dynamics \eqref{eq:RD}, the orbits of \eqref{eq:qRD} for $q>1$ may run into the boundary of the game's strategy space in finite time.
In this way, \eqref{eq:qRD} provides a smooth interpolation between the replicator dynamics and the projection dynamics (obtained for $q=1$ and $q=2$ respectively), with the replicator dynamics defining the boundary between the well- and ill-posed regimes of \eqref{eq:qRD}.
\end{example}

\begin{example}[The Rényi entropy]
The Rényi (negative) entropy is defined as
\begin{equation}
\label{eq:penalty-Renyi}
\txs
h(x)
	= - (1-q)^{-1} \log \left( \sum_{\alpha=1}^{n} x_{\alpha}^{q} \right)
\end{equation}
for $q \in (0,1)$.
Just like its Tsallis counterpart, the penalty function \eqref{eq:penalty-Renyi} is steep for all $q\in(0,1)$ and it approaches the Gibbs penalty function \eqref{eq:penalty-Gibbs} as $q\to1$.
Unlike \eqref{eq:penalty-Tsallis} though, \eqref{eq:penalty-Renyi} is not decomposable so we cannot use the explict formula \eqref{eq:RLD-theta} to derive the induced dynamics.
Still, after a somewhat tedious calculation (which we carry out in Appendix \ref{app:Renyi}), we obtain the \emph{Rényi dynamics}
\begin{equation}
\label{eq:ReD}
\tag{ReD}
\dot x_{k\alpha}
	= \frac{x_{k\alpha}}{\xi_{k\alpha}} \payv_{k\alpha}
	+ x_{k\alpha} \frac{S_{k,q} - \xi_{k\alpha}^{-1}}{1 - S_{k,q}} \insum_{\beta}^{k} x_{k\beta} \payv_{k\beta}
	- x_{k\alpha} \frac{1 - \xi_{k\alpha}^{-1}}{1- S_{k,q}} \insum_{\beta}^{k} \frac{x_{k\beta}}{\xi_{k\beta}} \payv_{k\beta},
\end{equation}
where
$\xi_{k\alpha} = q x_{k\alpha}^{q-1}\big/ \insum_{\gamma}^{k} x_{k\gamma}^{q}$
and
$S_{k,q} = q^{-1} \insum_{\gamma}^{k} x_{\gamma}^{q} \cdot \insum_{\gamma}^{k} x_{\gamma}^{2-q}$.

In view of its rather complicated form, it is important to recall that the system \eqref{eq:ReD} simply describes the evolution of the reinforcement learning dynamics \eqref{eq:RL} with the above choice of penalty function.
Furthermore, just as the Gibbs penalty \eqref{eq:penalty-Gibbs} is recovered from \eqref{eq:penalty-Renyi} in the limit $q\to1$, it is natural to expect that the replicator dynamics \eqref{eq:RD} may themselves be seen as a limiting case of \eqref{eq:ReD} as $q\to1$;
in Appendix \ref{app:Renyi} we show that this indeed the case.
\end{example}

\begin{example}[The log-barrier]
\label{ex:LD}
An important nonexample of a penalty function is the logarithmic barrier
\begin{equation}
\label{eq:penalty-log}
h(x)
	= -\sum_{\alpha=1}^{n} \log x_{\alpha}.
\end{equation}
Obviously, \eqref{eq:penalty-log} is steep, strongly convex and decomposable, but it is not finite at the boundary $\bd(\simplex)$ of $\simplex$.
Nevertheless,
letting $\theta(x) = -\log x$ (so $\theta''(x) = 1/x^{2}$) and working as in Example \ref{ex:qRD}, \eqref{eq:RLD-theta} yields the \emph{log-barrier dynamics}
\begin{equation}
\label{eq:LD}
\tag{LD}
\dot x_{k\alpha}
	= x_{k\alpha}^{2}\left[\payv_{k\alpha} - r_{k}^{-2}\insum_{\beta}^{k} x_{k\beta}^{2} \payv_{k\beta}\right],
\end{equation}
with $r_{k}^{2} = \insum_{\beta}^{k} x_{k\beta}^{2}$.
The system \eqref{eq:LD} is easily seen to be well-posed and it can be seen as a limiting case of \eqref{eq:qRD} when $q\to0^{+}$.
In convex optimization, \eqref{eq:LD} was first considered by \citeor{BL89} and it has since been studied extensively by many authors \textendash\ see e.g. \citeor{Fia90}, \citeor{Kiw97b}, \citeor{BT03}, \citeor{ABB04}, \citeor{LM15} and references therein.
The results that we derive in the rest of the paper for \eqref{eq:RL} remain true in the case of \eqref{eq:LD}, 
but we do not provide proofs.
\end{example}

\begin{figure}[t]
\subfigure{
\label{fig:portraits-PD}
\includegraphics[width=.48\textwidth]{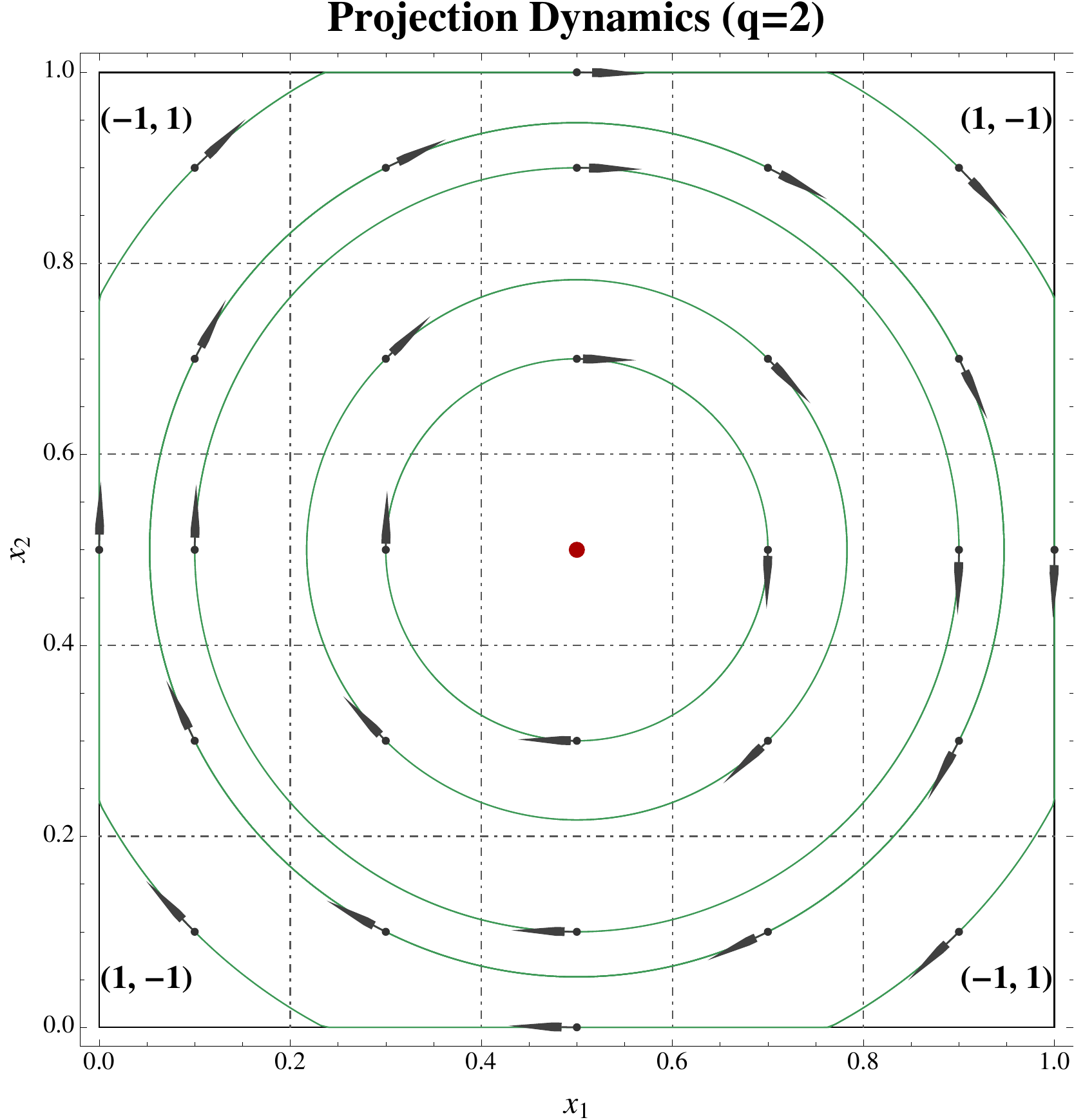}}
\hfill
\subfigure{
\label{fig:portraits-qRD}
\includegraphics[width=.48\textwidth]{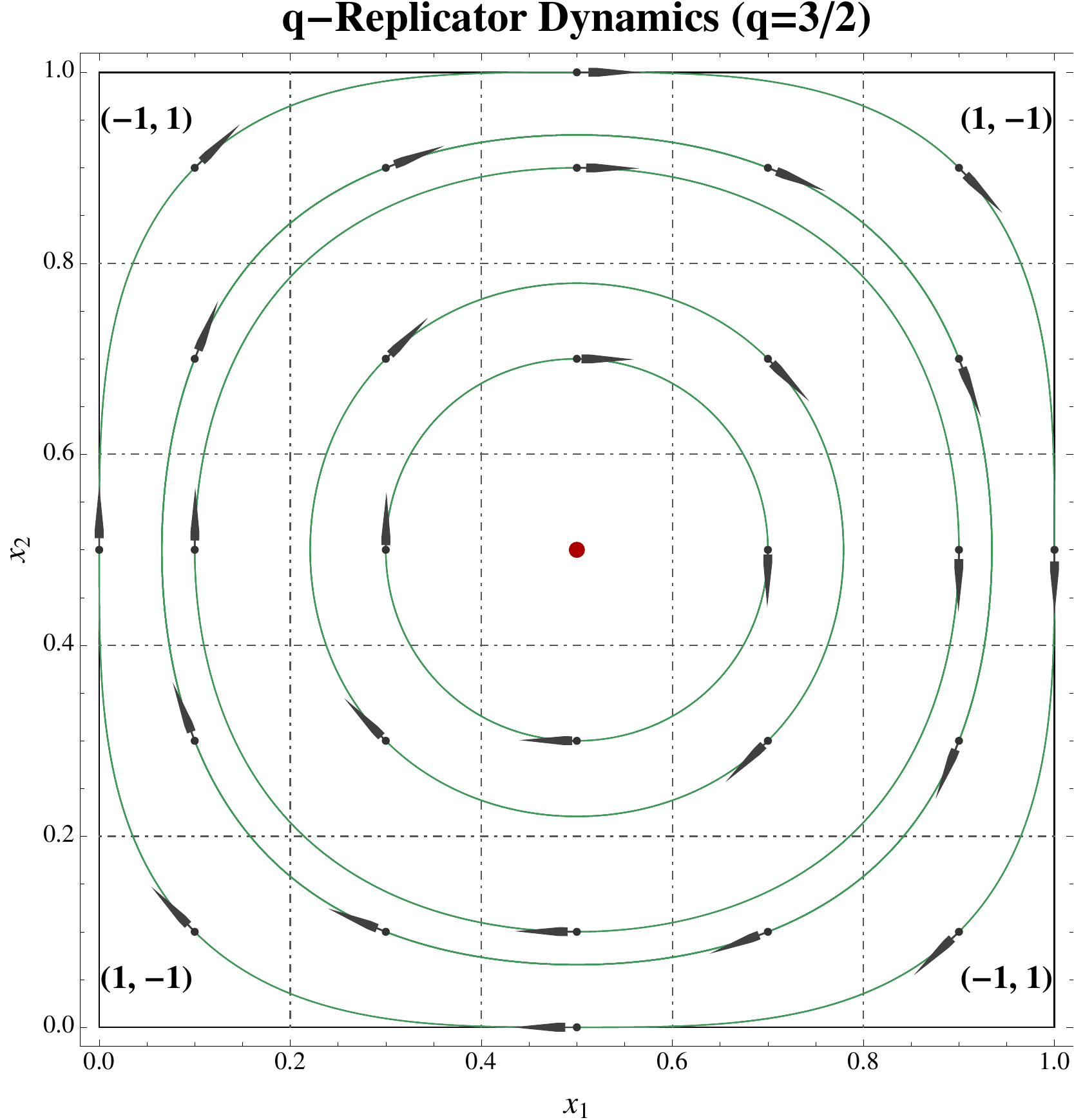}}
\quad

\subfigure{
\label{fig:portraits-RD}
\includegraphics[width=.48\textwidth]{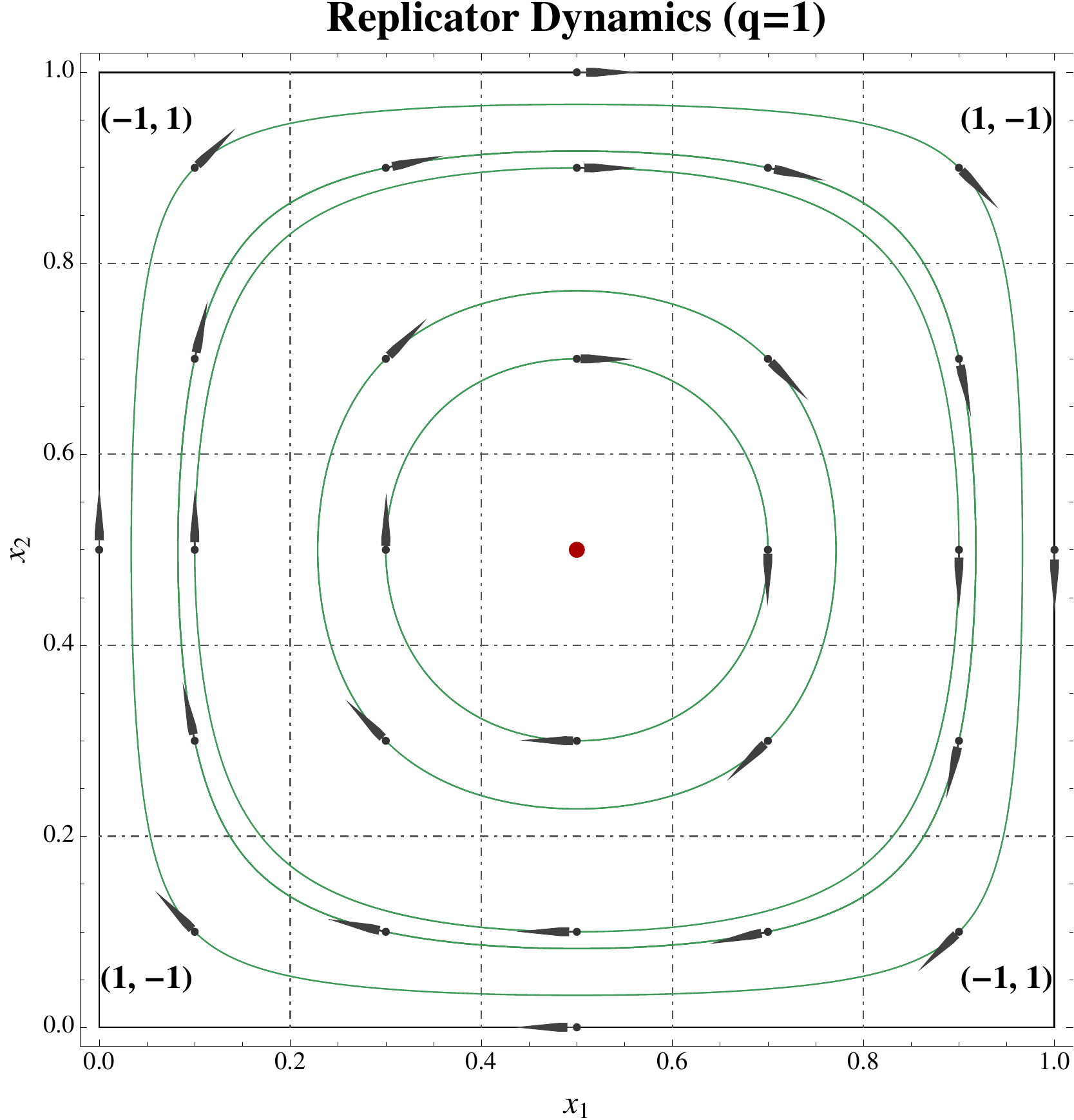}}
\hfill
\subfigure{
\label{fig:portraits-LD}
\includegraphics[width=.48\textwidth]{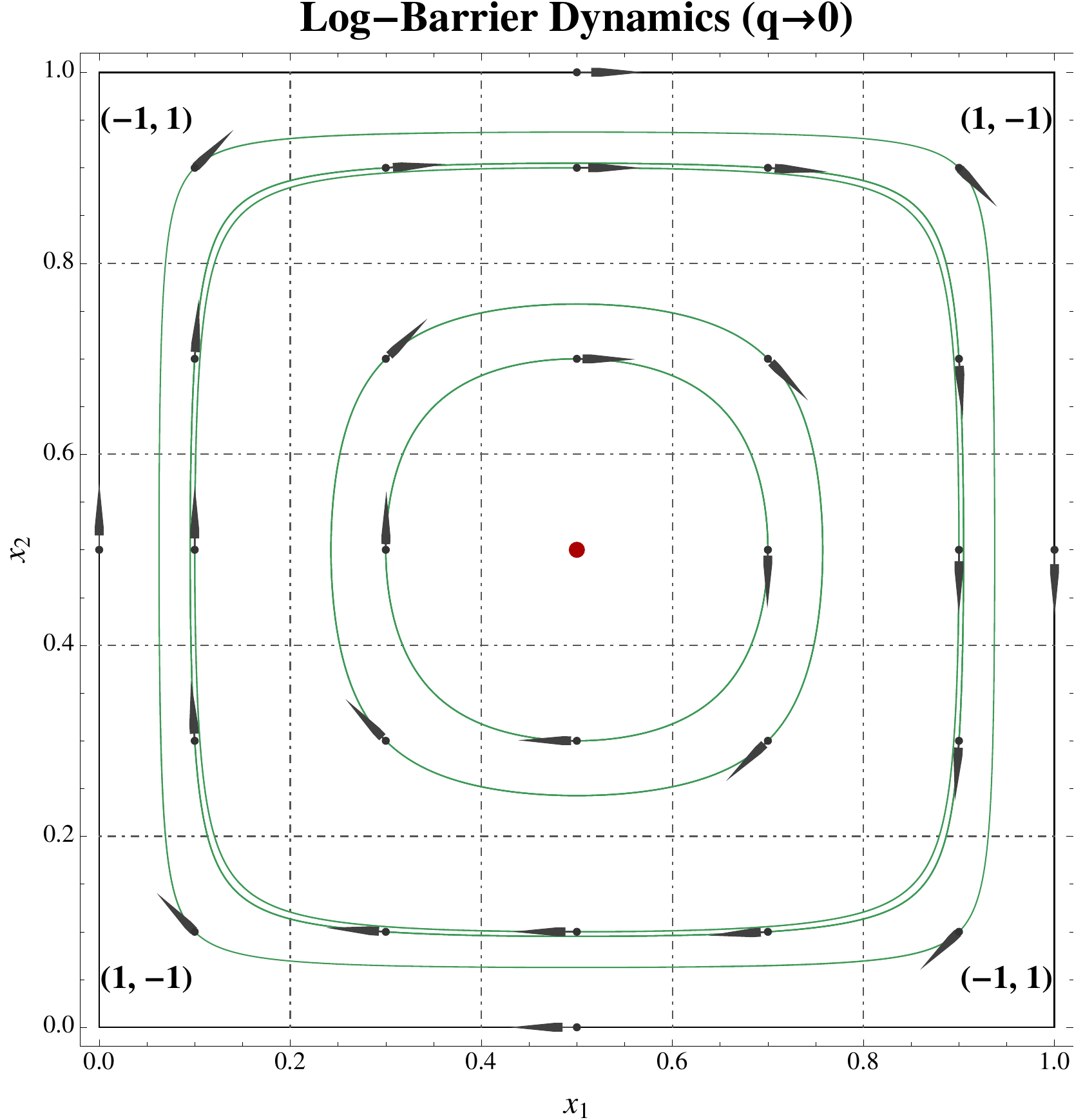}}
\quad

\caption{\small
Evolution of play under the reinforcement learning dynamics \eqref{eq:RL} in \acl{MP} (Nash equilibria are depicted in dark red and stationary points in light/dark blue; for the game's payoffs, see the vertex labels).
As the deformation parameter $q$ of \eqref{eq:qRD} decreases, we pass from the nonsteep regime ($q>1$) where the orbits of \eqref{eq:RL} in $\strat$ collide with the boundary of $\strat$ in finite time, to the steep regime ($q\leq1$) where \eqref{eq:qRD} becomes well-posed.}
\label{fig:portraits}
\end{figure}

\section{Elimination of dominated strategies}
\label{sec:dominance}

We begin our rationality analysis with the elimination of dominated strategies.
Formally, if $\game\equiv\game(\play,\act,\pay)$ is a finite game in normal form, we say that $p_{k}\in\strat_{k}$ is \emph{dominated} by $p_{k}'\in\strat_{k}$ and we write $p_{k}\prec p_{k}'$ when
\begin{equation}
\label{eq:dom-mixed}
\pay_{k}(p_{k};x_{-k}) < \pay_{k}(p_{k}';x_{-k})
	\quad
	\text{for all $x_{-k}\in\strat_{-k}\equiv\inprod_{\ell\neq k} \strat_{\ell}$.}
\end{equation}
Thus, for pure strategies $\alpha,\beta\in\act_{k}$, we have $\alpha\prec\beta$ whenever
\begin{equation}
\label{eq:dom-pure}
\payv_{k\alpha}(x)
	< \payv_{k\beta}(x)
	\quad
	\text{for all $x\in\strat$.}
\end{equation}
If \eqref{eq:dom-mixed} is strict for only some (but not all) $x\in\strat$, we will say that $p_{k}$ is \emph{weakly dominated} by $p_{k}'$ and we will write $p_{k} \preccurlyeq p_{k}'$.
Conversely, we will say that $q=(p_{1},\dotsc,p_{N})\in\strat$ is \emph{undominated} if no component $p_{k}\in\strat_{k}$ of $p$ is (strictly) dominated.
Of course, if dominated strategies strategies are removed from $\game$, other strategies may become dominated in the resulting restriction of $\game$, leading to the notion of \emph{iteratively dominated strategies}.
Accordingly, a strategy which survives all rounds of elimination is called \emph{iteratively undominated}.

For a given trajectory of play $x(t)\in\strat$, $t\geq0$, we say that the pure strategy $\alpha\in\act_{k}$ \emph{becomes extinct along $x(t)$} if $x_{k\alpha}(t)\to 0$ as $t\to\infty$.
More generally, following \citeor{SZ92}, we say that the mixed strategy $p_{k}\in\strat_{k}$ becomes extinct along $x(t)$ if $\min\{x_{k\alpha}(t): \alpha\in\supp(p_{k})\}\to0$;
otherwise, we say that $p_{k}$ \emph{survives}.

Extending the classic elimination results of \citeor{Aki80}, \citeor{Nac90}, and \citeor{SZ92}, we first show that only iteratively undominated strategies survive under \eqref{eq:RL}:

\begin{theorem}
\label{thm:dom}
Let $x(t) = \choice(y(t))$ be an orbit of \eqref{eq:RL} in $\strat$.
If $p_{k}\in\strat_{k}$ is dominated \textup(even iteratively\textup), then it becomes extinct along $x(t)$.
\end{theorem}

In the replicator dynamics, most proofs of elimination of dominated strategies involve some form of the \acl{KL} divergence function $\dkl(p_{k},x_{k}) = \insum_{\alpha}^{k} p_{k\alpha} \log (p_{k\alpha}/x_{k\alpha})$, an asymmetric measure of the ``distance'' between $p_{k}$ and $x_{k}$.
In particular, to show that $p_{k}$ is eliminated along $x_{k}(t)$ it suffices to show that $\dkl(p_{k},x_{k}(t))\to + \infty$.
Following \citeor{Bre67}, the same role for a steep penalty function $h\from\simplex\to\R$ is played by the so-called \emph{Bregman divergence}
\begin{equation}
\label{eq:Bregman-steep}
\breg_{h}(p, x)
	= h(p) - h(x) - \braket{dh(x)}{p-x},
	\qquad
	p\in\simplex,\:
	x\in\intsimplex,
\end{equation}
where $dh(x)$ denotes the differential of $h$ at $x$ (so $\breg_{h}(p,x)$ is just the difference between $h(p)$ and the estimate of $h(p)$ based on linearization at $x$).%
\footnote{One can easily verify that the Bregman divergence \eqref{eq:Bregman-steep} of the Gibbs penalty \eqref{eq:penalty-Gibbs} is simply the standard \acl{KL} divergence.}

On the other hand, since \eqref{eq:RLD} may fail to be well-posed if the players' penalty functions are not steep, we must analyze the reinforcement learning dynamics \eqref{eq:RL} directly on the dual space $\dual$ where the score variables $y$ evolve.
 To do so, we introduce here the \emph{Fenchel coupling} between $p_{k}\in \strat_{k}$ and $y_{k}\in \dual_{k}$, defined as
\begin{equation}
\label{eq:Fenchel}
\fench_{k}(p_{k},y_{k})
	= h_{k}(p_{k}) + h_{k}^{\ast}(y_{k}) - \braket{y_{k}}{p_{k}},
\end{equation}
where
\begin{equation}
h_{k}^{\ast}(y_{k})
	= \max_{x_{k}\in\strat_{k}}\{\braket{y_{k}}{x_{k}} - h_{k}(x_{k})\}
\end{equation}
denotes the convex conjugate of the penalty function $h_{k}\from\strat_{k}\to\R$ of player $k$.

Our choice of terminology above simply reflects the fact that $\fench_{k}(p_{k},y_{k})$ collects all the terms of Fenchel's inequality, so it is nonnegative and (strictly) convex in both arguments.
Furthermore, we show in Proposition \ref{prop:Bregman-dual} that
\begin{inparaenum}
[\itshape a\upshape)]
\item
$\fench_{k}(p_{k},y_{k})$ is equal to the associated Bregman divergence between $p_{k}$ and $x_{k} = \choice_{k}(y_{k})$ when the latter is interior;
and
\item
$\fench_{k}(p_{k},y_{k})$ provides a proximity measure between $p_{k}$ and $\choice_{k}(y_{k})$ which is applicable even when $h_{k}$ is not steep.
\end{inparaenum}

\begin{Proof}[Proof of Theorem \ref{thm:dom}]
Assume first that $p_{k}\in\strat_{k}$ is dominated by $p_{k}'\in\strat_{k}$ and let $\Lambda_{k} = \{x_{k}\in\strat_{k}: \text{$x_{k\alpha}=0$ for some $\alpha\in\supp(p_{k})$}\}$ be the union of all faces of $\strat_{k}$ that do not contain $p_{k}$.
By definition, $p_{k}$ becomes extinct along $x(t)$ if and only if $x_{k}(t) \to \Lambda_{k}$ as $t\to\infty$;
therefore, in view of Proposition \ref{prop:divergence}, it suffices to show that $\fench_{k}(p_{k},y_{k}(t))\to+\infty$.

To that end, consider the ``cross-coupling''
\begin{equation}
\label{eq:cross}
V_{k}(y_{k})
	= \fench_{k}(p_{k},y_{k}) - \fench_{k}(p_{k}',y_{k})
	= h_{k}(p_{k}) - h_{k}(p_{k}') - \braket{y_{k}}{p_{k} - p_{k}'}.
\end{equation}
Under the dynamics \eqref{eq:RL}, we then have:
\begin{equation}
\label{eq:cross1}
\frac{d}{dt} V_{k}(y_{k}(t))
	= -\braket{\payv_{k}(x(t))}{p_{k} - p_{k}'}
	= \pay_{k}(p_{k}';x_{-k}(t)) - \pay_{k}(p_{k};x_{-k}(t))
	\geq \delta_{k}
	> 0,
\end{equation}
where $\delta_{k} = \min_{x\in\strat}\{\pay_{k}(p_{k}';x_{-k}) - \pay_{k}(p_{k};x_{-k})\}$ denotes the minimum payoff difference between $p_{k}$ and $p_{k}'$.
Hence, with $\fench_{k}(p_{k}',y_{k})\geq0$ for all $y_{k}\in \dual_{k}$ (Proposition \ref{prop:Bregman-dual}), we readily obtain
\begin{equation}
\fench_{k}(p_{k},y_{k}(t))
	\geq V_{k}(y_{k}(0)) + \delta_{k} t,
\end{equation}
so every $\omega$-limit of $x_{k}(t)$ belongs to $\Lambda_{k}$ by Proposition \ref{prop:divergence}, i.e. $p_{k}$ becomes extinct.

To show that iteratively dominated strategies become extinct, we proceed by induction on the rounds of elimination of dominated strategies.
More precisely, let $\strat_{k}^{r}\subseteq\strat_{k}$ denote the space of mixed strategies of player $k$ that survive $r$ rounds of elimination so that all strategies $p_{k}\notin\strat_{k}^{r}$ become extinct along $x(t)$;
in particular, if $\alpha\notin\act_{k}^{r}\equiv\act_{k}\cap\strat_{k}^{r}$, this implies that $x_{k\alpha}(t)\to 0$ as $t\to\infty$.
Assume further that $p_{k}\in\strat_{k}^{r}$ survives for $r$ elimination rounds but dies on the subsequent one, so there exists some $p_{k}'\in\strat_{k}^{r}$ with $\pay_{k}(p_{k}';x_{-k}) > \pay_{k}(p_{k};x_{-k})$ for all $x\in\strat^{r} = \prod_{\ell} \strat_{\ell}^{r}$.
With this in mind, decompose $x\in\strat$ as $x=x^{r}+z^{r}$ where $x^{r}$ is the (Euclidean) projection of $x$ on the subspace of $\strat$ spanned by the surviving pure strategies $\act_{\ell}^{r}$, $\ell\in\play$.
Our induction hypothesis implies $z^{r}(t) = x(t) - x^{r}(t) \to 0$ as $t\to\infty$ (recall that $x_{k\alpha}(t)\to0$ for all $\alpha\notin\act_{k}^{r}$), so, for large enough $t$, we have
\begin{equation}
\label{eq:paydiff-dom}
\abs{%
	\braket{\payv_{k}(x(t))}{p_{k}'-p_{k}}
	- \braket{\payv_{k}(x^{r}(t))}{p_{k}'-p_{k}}
	}
	< \delta_{k}^{r}/2,
\end{equation}
where $\delta_{k}^{r} = \min_{x^{r}\in\strat^{r}} \{ \pay_{k}(p_{k}';x_{-k}^{r}) - \pay_{k}(p_{k};x_{-k}^{r}) \}$.

By combining the above, we get
\begin{equation}
\braket{\payv_{k}(x(t))}{p_{k}' - p_{k}}
	= \pay_{k}(p_{k}';x(t)) - \pay_{k}(p_{k};x(t))
	> \delta_{k}^{r}/2
	> 0
\end{equation}
for large $t$, and our claim follows by plugging this last estimate into \eqref{eq:cross1} and arguing as in the base case $r=1$.
\end{Proof}

\begin{remark}
In the projection dynamics of \citeor{NZ97}, dominated strategies need not be eliminated.
Although such strategies are selected against at interior states, the dynamics' solution orbits may enter and leave the boundary of $\strat$ in perpetuity, allowing dominated strategies to survive \citeorp{SDL08}.
Similarly, there exist Carath\'eodory solutions of the projection dynamics \eqref{eq:PD} that do not eliminate dominated strategies (for instance, stationary trajectories at vertices corresponding to dominated strategies).
By contrast, Theorem \ref{thm:dom} shows that dominated strategies become extinct along \emph{every} orbit $x(t)=\proj_{\strat} y(t)$ of \eqref{eq:PD} that is induced by the projected reinforcement learning scheme \eqref{eq:PL}.%
\footnote{Recall here that such orbits are solutions of \eqref{eq:PD} in the sense of Corollary \ref{cor:RLD-steep} \textendash\ i.e. they satisfy \eqref{eq:PD} over an open dense subset of $[0,\infty)$.}
\end{remark}

\begin{remark}
Theorem \ref{thm:dom} imposes no restrictions on the players' choice maps.
For instance, dominated strategies become extinct even if some players use the exponential learning scheme \eqref{eq:XL} while others employ the projection-driven process \eqref{eq:PL}.  
In fact, the proof of Theorem \ref{thm:dom} shows that the elimination of a player's dominated strategies is a unilateral result:
if a player follows \eqref{eq:RL}, he ceases to play dominated strategies irrespective of what other players are doing.
\end{remark}

We now turn to the \emph{rate} of elimination of dominated strategies.
In the case of the replicator dynamics, this rate is known to be exponential: if $\alpha\prec\beta$, then $x_{\alpha}(t) = \bigoh(\exp(-ct))$ for some $c>0$ \textendash\ see e.g.~\citeor{Wei95}.
However, as we show below, the rate of elimination of dominated strategies under \eqref{eq:RL} depends crucially on the players' choice of penalty function;
in fact, if the players' penalty functions are nowhere steep, dominated strategies become extinct in \emph{finite} time.
This is perhaps most easily seen in the case of the projection dynamics \eqref{eq:PL}:
when starting at the interior of the simplex, a player simply follows his payoff vector and there is no mitigating factor to slow down his trajectory of play near the boundary of the simplex (by contrast, in the replicator dynamics, $\dot x_{k\alpha}$ tends to zero as $x_{k\alpha}\to0$).
In this way, trajectories of play attain the boundary of the simplex in finite time \textendash\ and so on for every subface of the simplex until all dominated strategies are eliminated.

Building on this intuition, our general result is as follows:

\begin{proposition}
\label{prop:dom-rate}
Let $x(t) = \choice(y(t))$ be an orbit of the dynamics \eqref{eq:RL-rate} and assume that the players' penalty functions are of the form $h_{k}(x_{k}) = \insum_{\beta}^{k} \theta_{k}(x_{k\beta})$ for some $\theta_{k}\from[0,1]\to\R$ as in \eqref{eq:decomposable}.
If $\alpha\prec\beta$, then
\begin{equation}
\label{eq:dom-rate}
x_{k\alpha}(t)
	\leq \rate_{k}(c_{k} - \tempinv_{k} \delta_{k}t),
\end{equation}
where
$c_{k}$ is a constant that only depends on the initial conditions of \eqref{eq:RL},
$\delta_{k} = \min\{\payv_{k\beta}(x) - \payv_{k\alpha}(x):x\in\strat\}$ is the minimum payoff difference between $\alpha$ and $\beta$,
and the rate function $\rate_{k}$ is given by:
\begin{equation}
\label{eq:rate}
\rate_{k}(z)
	= \begin{cases}
	0
		&\quad
		\text{if $z\leq \theta_{k}'(0^{+})$,}
	\\
	1
		&\quad
		\text{if $z\geq \theta_{k}'(1^{-})$,}
	\\
	(\theta_{k}')^{-1}(z)
		&\quad
		\text{otherwise,}
	\end{cases}
\end{equation}
where $(\theta_{k}')^{-1}$ is the inverse function of $\theta_{k}'$.
In particular, if $\theta_{k}'(0)$ is finite, dominated strategies become extinct in finite time.
\end{proposition}

\begin{Proof}
By the definition of the reinforcement learning dynamics \eqref{eq:RL-rate}, we have:
\begin{equation}
\label{eq:rate1}
\dot y_{k\alpha} - \dot y_{k\beta}
	= \tempinv_{k}\left[ \payv_{k\alpha}(x(t)) - \payv_{k\beta}(x(t)) \right]
	\leq - \tempinv_{k} \delta_{k},
\end{equation}
and hence:
\begin{equation}
\label{eq:rate2}
y_{k\alpha}(t) - y_{k\beta}(t)
	\leq y_{k\alpha}(0) - y_{k\beta}(0) - \tempinv_{k} \delta_{k} t.
\end{equation}
On the other hand, by the \ac{KKT} conditions \eqref{eq:KKT} for the softmax problem \eqref{eq:softmax}, we obtain $\theta_{k}'(x_{k\alpha}) - \theta_{k}'(x_{k\beta}) \leq y_{k\alpha} - y_{k\beta}$ whenever $x_{k\alpha} = \choice_{k\alpha}(y_{k}) > 0$.
Since $\theta_{k}'$ is bounded above on $(0,1]$, \eqref{eq:rate2} gives
\begin{equation}
\label{eq:rate3}
\theta_{k}'(x_{k\alpha}(t))
	\leq c_{k} - \tempinv_{k} \delta_{k}t
\end{equation}
for some $c_{k}\in\R$ and for all $t$ such that $x_{k\alpha}(t)>0$, so \eqref{eq:dom-rate} follows \textendash\ simply recall that $\theta_{k}'(0) \leq \theta_{k}'(x_{k\alpha}(t))$ so $x_{k\alpha}(t) = 0$ if $t \geq (c_{k} - \theta_{k}'(0))/(\tempinv_{k}\delta_{k})$.
\end{Proof}

\begin{remark}
Proposition \ref{prop:dom-rate} shows that a player's penalty function can be re\-verse-en\-gi\-neered in terms of the desired rate of elimination of dominated strategies:
to achieve a target extinction rate $\rate$, it suffices to pick a penalty kernel $\theta$ such that $\theta' = \rate^{-1}$ (cf. Table \ref{tab:rates}).
For instance, the Gibbs kernel $\theta(x) = x\log x$ of \eqref{eq:penalty-Gibbs} yields the exponential extinction rate $\exp(c-\delta t)$ whereas the quadratic kernel \eqref{eq:penalty-Eucl} gives the bound $x_{k\alpha}(t) \leq [c-\delta t]_{+}$ which shows that \eqref{eq:PL} eliminates dominated strategies in finite time.
\end{remark}

\begin{table}[tbp]
\centering
\renewcommand{\arraystretch}{1.4}
{\sc
\small
\begin{tabular}{llll}
\hline
\noalign{\vspace{1pt}}
\textsc{Dynamics}
	&\hspace{3em}
	&\textsc{Penalty Kernel $\theta(x)$}
	\hspace{1em}
	&\textsc{Decay Rate $\rate(-y)$}
	\\
\hline
\hline
projection
	&\eqref{eq:PL}
	&$\frac{1}{2}x^{2}$
	&$-y$
	\\
\hline
replicator
	&\eqref{eq:RD}
	&$x \log x$
	&$\exp(-y)$
	\\
\hline
{\small$q$}-replicator
	&\eqref{eq:qRD}
	&$[q(1-q)]^{-1} (x - x^{q})$
	&$[q^{-1} + (1-q) y]^{1/(q-1)}$
	\\
\hline
log-barrier
	&\eqref{eq:LD}
	&$-\log x$
	&$1/y$
	\\
\hline
\end{tabular}
}
\vspace{2ex}
\caption{\small
Rates of extinction of dominated strategies and convergence to strict equilibria under the dynamics \eqref{eq:RL-rate} for different penalty kernels $\theta$.
If $\alpha\prec\beta$ and $\payv_{k\beta} - \payv_{k\alpha} \geq \delta$, the rate function $\rate$ is such that $x_{k\alpha}(t) \leq \rate(c-\tempinv\delta t)$ for some $c\in\R$ \textendash\ cf. \eqref{eq:rate}.
Otherwise, if \revised{$\eq = (\alpha_{1}^{\ast},\dotsc,\alpha_{N}^{\ast})$ is a strict equilibrium and $\payv_{k,\alpha_{k}^{\ast}}(\eq) - \payv_{k\beta}(\eq) \geq \delta > 0$ for all $\beta\in\act_{k}\exclude\{\alpha_{k}^{\ast}\}$, the rate function $\rate$ is such that $1 - x_{k,\alpha_{k}^{\ast}}(t)$ is of the order of $\rate(c - \tempinv\delta t)$ for large $t$.}
If the penalty kernel $\theta$ is not steep at $0$ (as in the projection and $q$-replicator dynamics for $q>1$), extinction of dominated strategies and/or convergence to a strict equilibrium occurs in finite time.}
\label{tab:rates}
\end{table}

Finally, for weakly dominated strategies, we obtain the following conditional extinction result in the spirit of \citeor[Proposition 3.2]{Wei95}:

\begin{proposition}
\label{prop:dom-weak}
Let $x(t) = \choice(y(t))$ be an orbit of \eqref{eq:RL} in $\strat$ and let $p_{k}\preccurlyeq p_{k}'$.
Then, $p_{k}$ becomes extinct along $x(t)$ or every $\alpha_{-k}\in\act_{-k}\eqdef\prod_{\ell\neq k}\act_{\ell}$ such that $\pay_{k}(p_{k};\alpha_{-k}) < \pay_{k}(p_{k}';\alpha_{-k})$ becomes extinct along $x(t)$.
\end{proposition}

\begin{remark}
The ``or'' above is not exclusive:
the two extinction clauses could occur simultaneously \textendash\ see e.g. \citeor[Proposition 3.2]{Wei95} for the case of the replicator dynamics.
\end{remark}

\begin{Proof}[Proof of Proposition \ref{prop:dom-weak}]
With notation as in the proof of Theorem \ref{thm:dom}, we have:
\begin{equation}
\dot V_{k}
	= \braket{\payv_{k}(x(t))}{p_{k}' - p_{k}}
	= \insum_{\alpha_{-k}\in\act_{-k}'}
	\left[
	\pay_{k}(p_{k}';\alpha_{-k}) - \pay_{k}(p_{k};\alpha_{-k})
	\right]
	x_{\alpha_{-k}}(t),
\end{equation}
where $x_{\alpha_{-k}} \eqdef \prod_{\ell\neq k} x_{\ell,\alpha_{\ell}}$ denotes the $\alpha_{-k}$-th component of $x$ and $\act_{-k}' = \{\alpha_{-k}\in\act_{k}:\pay_{k}(p_{k};\alpha_{-k}) < \pay_{k}(p_{k}';\alpha_{-k})\}$.
Integrating with respect to $t$, we see that $V_{k}(y(t))$ remains bounded if and only if the integrals $\int_{0}^{\infty} x_{\alpha_{-k}}(t) \dd t$ are all finite.
However, with $\dot x_{\alpha_{-k}}(t)$ essentially bounded,
the same argument as in the proof of \citeor[Prop.~3.2]{Wei95} shows that $\lim_{t\to\infty} x_{\alpha_{-k}}(t) = 0$ if $\int_{0}^{\infty} x_{\alpha_{-k}}(t) \dd t$ is finite.
If this is not the case, we have $\lim_{t\to\infty} V_{k}(y(t)) = +\infty$ and Proposition \ref{prop:divergence} shows that $p_{k}$ becomes extinct.
\end{Proof}

\section{Equilibrium, stability and convergence}
\label{sec:folk}

We now turn to the long-term stability and convergence properties of the reinforcement dynamics \eqref{eq:RL}.
Our analysis focuses on \emph{Nash equilibria}, i.e. strategy profiles $\eq=(\eq_{1},\dotsc,\eq_{N})\in\strat$ that are unilaterally stable in the sense that
\begin{equation}
\label{eq:Nash}
\pay_{k}(\eq)
	\geq \pay_{k}(x_{k};\eq_{-k})
	\quad
	\text{for all $x_{k}\in\strat_{k}$ and for all $k\in\play$,}
\end{equation}
or, equivalently:
\begin{equation}
\label{eq:Nash-components}
\payv_{k\alpha}(\eq)
	\geq \payv_{k\beta}(\eq)
	\quad
	\text{for all $\alpha\in\supp(\eq_{k})$ and for all $\beta\in\act_{k}$, $k\in\play$.}
\end{equation}
If \eqref{eq:Nash} is strict for all $x_{k}\neq\eq_{k}$, $k\in\play$, we say that $\eq$ is a \emph{strict equilibrium}.
Finally, equilibria of restrictions of $\game$ are called \emph{restricted equilibria} of $\game$.
In particular, $\eq$ is a restricted equilibrium of $\game$ if \eqref{eq:Nash} holds for every $k\in\play$ and for all $x_{k}$ with $\supp(x_{k})\subseteq\supp(\eq_{k})$.

Some basic long-term stability and convergence properties of the replicator dynamics for (asymmetric) normal form games
can be summarized as follows:
\begin{enumerate}
\setlength{\itemsep}{1pt}
\setlength{\parskip}{1pt}
\item
Nash equilibria are stationary.
\item
If an interior solution orbit converges,
\revised{its limit is a Nash equilibrium.}
\item
If a point is Lyapunov stable,
\revised{then it is a Nash equilibrium.}
\item
\label{itm:strict}
Strict equilibria are asymptotically stable.
\end{enumerate}
Our aim in this section is to establish analogous results for the reinforcement learning scheme \eqref{eq:RL}.
That said, since \eqref{eq:RL} does not evolve directly on $\strat$ (and the induced dynamics \eqref{eq:RLD} are well-posed only when the players' penalty functions are steep), the standard notions of stability and stationarity must be modified accordingly.%
\footnote{The standard stability notions continue to apply in the dual space $\dual$ where $y(t)$ evolves;
however, since the mapping $\choice\from \dual\to\strat$ which defines the trajectories of play $x(t) = \choice(y(t))$ in $\strat$ is neither injective nor surjective, this approach would not suffice to define stationarity and stability on $\strat$.}

\begin{definition}
\label{def:stability}
Let $\eq\in\strat$ and let $x(t) = \choice(y(t))$ be a solution orbit of \eqref{eq:RL}.
We will say that:
\begin{enumerate}
\addtolength{\itemsep}{1pt}

\item
$\eq$ is \emph{stationary} under \eqref{eq:RL} if
$\eq \in\im\choice$
and $x(t) = \eq$ for all $t\geq0$ whenever $x(0) = \eq$.

\item
$\eq$ is \emph{Lyapunov stable} under \eqref{eq:RL} if, for every neighborhood $U$ of $\eq$, there exists a neighborhood $V$ of $\eq$ such that \revised{$x(t)\in U$} for all $t\geq0$ whenever $x(0)\in V\cap \im \choice$.

\item
$\eq$ is \emph{attracting} under \eqref{eq:RL} if it admits a neighborhood $V$ such that $x(t) \to \eq$ as $t\to\infty$ whenever $x(0)\in V\cap\im\choice$.

\item
$\eq$ is \emph{asymptotically stable} under \eqref{eq:RL} if it is Lyapunov stable and attracting.
\end{enumerate}
\end{definition}

\begin{remark}
The requirement $x(0)\in V\cap\im\choice$ above is redundant because $x(0)\in\im\choice$ by definition.
We only mention it to clarify that there are boundary points of $\strat$ which may be inadmissible as initial points of the dynamics \eqref{eq:RL}.

On a similar note, stationary points $\eq\in\strat$ are explicitly required to belong to the image of the players' choice map $\choice$ but no such assumption is made for stable states.
From a propositional point of view, this is done to ensure that points $\eq\in\strat\exclude\im\choice$ are not called stationary vacuously.
From a dynamical standpoint,
stationary points should themselves be (constant) trajectories of the dynamical system under study, whereas Lyapunov stable and attracting states only need to be \emph{approachable} by trajectories.

Since $\im\choice\supseteq\intstrat$, any point in $\strat$ can be a candidate for (asymptotic) stability under \eqref{eq:RL}.
However, boundary points might not be suitable candidates for stationarity, so stability does \emph{not} imply stationarity (as would be the case for a dynamical system defined on $\strat$).
In particular, recall that $\eq\notin\im\choice$ if and only if some player's penalty function is steep at $\eq$.
As such, in the (steep) example of exponential learning, $\eq\in\strat$ is stationary under \eqref{eq:XL} if and only if it is an \emph{interior} stationary point of the replicator dynamics \eqref{eq:RD}.
By contrast, in the (nonsteep) projection setting of \eqref{eq:PD}, \emph{any} point in $\strat$ may be stationary.
\end{remark}

With this definition at hand, we then obtain:

\begin{theorem}
\label{thm:folk}
Let $\game\equiv\game(\play,\act,\pay)$ be a finite game, let $\eq\in\strat$, and let $x(t) = \choice(y(t))$ be an orbit of \eqref{eq:RL} in $\strat$.
\begin{enumerate}
[\upshape I.]
\addtolength{\itemsep}{1pt}

\item
\label{itm:folk1}
If $\eq$ is stationary under \eqref{eq:RL}, then it is a Nash equilibrium of $\game$;
conversely, if $\eq$ is a Nash equilibrium of $\game$ and $\eq\in\im\choice$, $\eq$ is stationary under \eqref{eq:RL}.

\item
\label{itm:folk2}
If\, $\lim_{t\to\infty} x(t) = \eq$, then $\eq$ is a Nash equilibrium of $\game$.

\item
\label{itm:folk3}
If $\eq\in\strat$ is Lyapunov stable under \eqref{eq:RL}, then $\eq$ is a Nash equilibrium of $\game$.

\item
\label{itm:folk4}
If $\eq$ is a strict Nash equilibrium of $\game$, then it is also asymptotically stable under \eqref{eq:RL}.
\end{enumerate}
\end{theorem}

%

For steep and decomposable penalty functions, Parts \ref{itm:folk1}, \ref{itm:folk3} and \ref{itm:folk4} of Theorem \ref{thm:folk} essentially follow from Theorem 1 in \citeor{CGM15} (see also \citeauthor{LM13} \cite{LM13,LM15} for related results in a second order setting).
Our proofs mimic those of \citeor{CGM15}, but the lack of steepness means that we must work directly on the dual space $\dual$ of the score variables $y_{k}$ and rely on the properties of the Fenchel coupling.

To prove Theorem \ref{thm:folk}, we need the following result (which is of independent interest):

\begin{proposition}
\label{prop:stability}
If every neighborhood $U$ of $\eq\in\strat$ admits an orbit $x_{U}(t) = \choice(y_{U}(t))$ of \eqref{eq:RL} such that $x_{U}(t) \in U$ for all $t\geq0$, then $\eq$ is a Nash equilibrium.
\end{proposition}

\begin{Proof}
Assume ad absurdum that $\eq$ is not \revised{a Nash equilibrium}, so $\payv_{k\alpha}(\eq) < \payv_{k\beta}(\eq)$ for some player $k\in\play$ and for some $\alpha\in\supp(\eq_{k})$, $\beta\in\act_{k}$.
Moreover, let $U$ be a sufficiently small neighborhood of $\eq$ in $\strat$ such that $\payv_{k\beta}(x) - \payv_{k\alpha}(x) \geq \delta$ for some $\delta > 0$ and for all $x\in U$.
Then, if $x(t) = \choice(y(t))$ is an orbit of \eqref{eq:RL} in $\strat$ that is contained in $U$ for all $t\geq0$, we get
\begin{equation}
y_{k\alpha}(t) - y_{k\beta}(t)
	= y_{k\alpha}(0) - y_{k\beta}(0)
	+ \int_{0}^{t} \left[\payv_{k\alpha}(x(s)) - \payv_{k\beta}(x(s))\right] \dd s
	\leq c - \delta t,
\end{equation}
for some $c\in\R$ and for all $t\geq0$.
This shows that $y_{k\alpha}(t) - y_{k\beta}(t)\to-\infty$ so, by Proposition \ref{prop:choice} (in Appendix \ref{app:choice}), we must also have $\lim_{t\to\infty} x_{k\alpha}(t) = 0$.
This contradicts the assumption that $x(t)$ remains in a small enough neighborhood of $\eq$ (recall that $\alpha\in\supp(\eq_{k})$ by assumption), so $\eq$ must be \revised{a Nash equilibrium}.
\end{Proof}

With this result at hand, we may proceed with the proof of Theorem \ref{thm:folk}:

\smallskip

\begin{Proof}[Proof of Theorem \ref{thm:folk}]
\hfill

\smallskip
\subparagraph{\sc Part \ref{itm:folk1}}
\quad
If $\eq$ is stationary under \eqref{eq:RL}, then $\choice(y(t)) = \eq$ for some $y(0)\in \dual$ and for all $t\geq0$;
this shows that $\eq$ satisfies the hypothesis of Proposition \ref{prop:stability}, so $\eq$ must be a Nash equilibrium of $\game$.
Conversely, assume that $\eq$ is \revised{a Nash equilibrium} and $\eq = \choice(y(0))$ for some initial $y(0)\in \dual$;
we then claim that the trajectory $y(t)$ with $y_{k\alpha}(t) = y_{k\alpha}(0) + \payv_{k\alpha}(\eq) t$ is the unique solution of \eqref{eq:RL} starting at $y(0)$.
Indeed, since $\payv_{k\beta}(\eq) \leq \payv_{k\alpha}(\eq)$ for all $\alpha\in\supp(\eq_{k})$ and for all $\beta\in\act_{k}$, we have $y_{k\alpha}(t) = y_{k\alpha}(0) + c_{k}t - d_{k\alpha}t$ where $d_{k\alpha}=0$ if $\alpha\in\supp(\eq_{k})$ and $d_{k\alpha}\geq0$ otherwise.
Proposition \ref{prop:choice} shows that $\choice_{k}(y_{k}(t)) = \eq_{k}$, so $y(t)$ satisfies \eqref{eq:RL} and our assertion follows by the well-posedness of \eqref{eq:RL}.

\medskip
\subparagraph{\sc Parts \ref{itm:folk2} and \ref{itm:folk3}}
\quad
The convergence and stability assumptions of Parts \ref{itm:folk2} and \ref{itm:folk3} both imply the hypothesis of Proposition \ref{prop:stability} so $\eq$ must be a Nash equilibrium of $\game$.

\medskip
\subparagraph{\sc  Part \ref{itm:folk4}}
\quad
Let $\eq = (\alpha_{1}^{\ast},\dotsc,\alpha_{N}^{\ast})$ be a strict equilibrium of $\game$, let $\act_{k}^{ \ast} = \act_{k}\exclude\{\alpha_{k}^{\ast}\}$, and consider the relative score variables
\begin{equation}
\label{eq:zscore}
z_{k\mu}
	= y_{k\mu} - y_{k\alpha_{k}^{\ast}},
	\quad
	\mu\in\act_{k}^{\ast},
\end{equation}
so that
\begin{equation}
\label{eq:RL-Z}
\dot z_{k\mu}
	= \payv_{k\mu}(x) - \payv_{k\alpha_{k}^{\ast}}(x).
\end{equation}
Proposition \ref{prop:choice} shows that $x_{k\mu}\to0$ whenever $z_{k\mu}\to-\infty$, so we also have $x\to\eq$ if $z_{k\mu}\to-\infty$ for all $\mu\in\act_{k}^{\ast}$, $k\in\play$.
Moreover, given that $\eq$ is a strict equilibrium, the RHS of \eqref{eq:RL-Z} is negative if $x$ is close to $\eq$;
the main idea of our proof will thus be to show that the relative scores $z_{k\mu}$ escape to negative infinity when they are not too large to begin with (i.e. when $x(0)$ is close enough to $\eq$).

To make this precise, let $\delta_{k} = \min_{\mu\in\act_{k}^{\ast}}\{\payv_{k\alpha_{k}^{\ast}}(\eq) - \payv_{k\mu}(\eq)\} > 0$ and let $\eps>0$ be such that $\payv_{k\mu}(x) - \payv_{k\alpha_{k}^{\ast}}(x) < - \delta_{k}/2$ for all $x$ with $\norm{x - \eq}^{2} < \eps$.
Furthermore, let
\begin{equation}
\fench_{h}(\eq,y)
	= \insum_{k} \fench_{k}(\eq_{k},y_{k})
	= \insum_{k} \left[ h_{k}(\eq_{k}) + h_{k}^{\ast}(y_{k}) - \braket{y_{k}}{\eq_{k}}\right]
\end{equation}
denote the Fenchel coupling between $\eq$ and $y$ (cf. Appendix \ref{app:Bregman}),
and set
\begin{equation}
\label{eq:nhd}
U_{\eps}^{\ast}
	= \left\{y\in \dual: \fench_{h}(\eq,y) < \eps K_{\min}/2\right\}
\end{equation}
where $K_{\min}>0$ is the smallest strong convexity constant of the players' penalty functions.
Proposition \ref{prop:Bregman-dual} gives $\norm{\choice(y) - \eq}^{2} \leq 2 K_{\min}^{-1} \fench_{h}(\eq,y) < \eps$, so we also have $\payv_{k\mu}(\choice(y)) - \payv_{k\alpha_{k}^{\ast}}(\choice(y)) < - \delta_{k}/2$ for all $y\in U_{\eps}^{\ast}$.

In view of the above, let $y(t)$ be a solution of \eqref{eq:RL} with $y(0)\in U_{\eps}^{\ast}$ and let $\tau_{\eps} = \inf\{t:y(t) \notin U_{\eps}^{\ast}\}$ be the first exit time of $y(t)$ from $U_{\eps}^{\ast}$.
Then, if $\tau_{\eps}<\infty$:
\begin{equation}
\label{eq:zexit}
z_{k\mu}(\tau_{\eps})
	= z_{k\mu}(0) + \int_{0}^{\tau_{\eps}} \left[ \payv_{k\mu}(x(s)) - \payv_{k\alpha_{k}^{\ast}}(x(s)) \right] \dd s
	\leq z_{k\mu}(0) - \tfrac{1}{2} \delta_{k} \tau_{\eps}
	< z_{k\mu}(0),
\end{equation}
for all $\mu\in\act_{k}^{\ast}$ and for all $k\in\play$.
Intuitively, since the score differences $y_{k\alpha_{k}^{\ast}} - y_{k\mu}$ grow with $t$, $x(\tau_{\eps})$ must be closer to $\eq$ than $x(0)$, meaning that $y(\tau_{\eps}) \in U_{\eps}^{\ast}$, a contradiction.
More rigorously, note that $\fench_{k}(\eq_{k},y_{k}) = h_{k}(\eq_{k}) + h_{k}^{\ast}(y_{k}) - y_{k\alpha_{k}^{\ast}}$.
Since $\braket{y_{k}}{x_{k}} = y_{k\alpha_{k}^{\ast}} + \sum_{\mu\in\act_{k}^{\ast}} x_{k\mu} (y_{k\mu} - y_{k\alpha_{k}^{\ast}})$, we will also have $h_{k}^{\ast}(y_{k}) - y_{k\alpha_{k}^{\ast}} = \max_{x_{k}\in\strat_{k}} \{ \insum_{\mu\in\act_{k}^{\ast}} x_{k\mu} z_{k\mu} - h_{k}(x_{k}) \}$, so $h_{k}^{\ast}(y_{k}') - y_{k\alpha_{k}^{\ast}}' < h_{k}^{\ast}(y_{k}) - y_{k\alpha_{k}^{\ast}}$ whenever $z_{k\mu}' < z_{k\mu}$ for all $\mu\in\act_{k}^{\ast}$.
In this way, \eqref{eq:zexit} yields the contradictory statement $\eps \leq \fench_{h}(q,y(\tau_{\eps})) < \fench_{h}(q,y(0)) < \eps$, so $y(t)$ must remain in $U_{\eps}^{\ast}$ for all $t\geq0$.
The estimate \eqref{eq:zexit} then shows that $\lim_{t\to\infty} z_{k\mu}(t) = -\infty$, so $\lim_{t\to\infty} x_{k\mu}(t) = 0$ by Proposition \ref{prop:choice}.
Therefore, $\lim_{t\to\infty} x(t) = \eq$.

The above shows that $y(t)$ is contained in $U_{\eps}^{\ast}$ and $\choice(y(t)) \to \eq$ whenever $y(0) \in U_{\eps}^{\ast}$.
To complete the proof, let $U_{\eps} = \{x\in\strat: \insum_{k} \breg_{h_{k}}(\eq_{k},x_{k}) < \eps K_{\min}/2\}$ where $\breg_{h_{k}}(\eq_{k},x_{k})$ is the Bregman divergence \eqref{eq:Bregman} between $x_{k}$ and $\eq_{k}$.
Propositions \ref{prop:Bregman} and \ref{prop:Bregman-dual} show that $U_{\eps}$ is a neighborhood of $\eq$ in $\strat$ with $\choice(U_{\eps}^{\ast}) \subseteq U_{\eps}$ and $\choice^{-1}(U_{\eps}) = U_{\eps}^{\ast}$, 
so $\eq$ is asymptotically stable under \eqref{eq:RL}.
%
\end{Proof}

\begin{remark}
In the (asymmetric) replicator dynamics \eqref{eq:RD}, it is well known that only strict equilibria can be attracting \textendash\ hence strict Nash equilibria and asymptotically stable states coincide.
One can extend this equivalence to \eqref{eq:RL} by using restrictions of $\game$ with smaller strategy sets to define the mixed strategy dynamics \eqref{eq:RLD} on the faces of $\strat$ \textendash\ for a related discussion, see \citeor{CGM15}.
\end{remark}

Theorem \ref{thm:folk} shows that the reinforcement learning scheme \eqref{eq:RL} exhibits essentially the same long-run properties as the benchmark replicator dynamics.
That said, from a quantitative viewpoint, the situation can be quite different:
as we show below, the rate of convergence of \eqref{eq:RL} to strict equilibria depends crucially on the players' penalty functions, and convergence can occur in finite time.
More formally, we have:

\begin{proposition}
\label{prop:strict-rate}
Let $x(t) = \choice(y(t))$ be an orbit of \eqref{eq:RL-rate},
let $\eq = (\alpha_{1}^{\ast},\dotsc,\alpha_{N}^{\ast})$ be a strict Nash equilibrium,
and assume that the players' penalty functions are of the form $h_{k}(x_{k}) = \insum_{\beta}^{k} \theta_{k}(x_{k\beta})$ for some $\theta_{k}\from[0,1]\to\R$ as in \eqref{eq:decomposable}.
\revised{%
Then, for every $\eps>0$ and for all $x(0)$ sufficiently close to $\eq$, we have:
\begin{equation}
\label{eq:conv-rate}
1 - x_{k\alpha_{k}^{\ast}}(t)
	\leq \sum_{\mu\in\act_{k}\exclude\{\alpha_{k}^{\ast}\}}
	\rate_{k}(c_{k\mu} - (1-\eps) \tempinv_{k}\delta_{k\mu} t),
\end{equation}
where
$c_{k\mu}$ is a constant that only depends on the initial conditions of \eqref{eq:RL},
$\delta_{k\mu} = \payv_{k\alpha_{k}^{\ast}}(\eq) - \payv_{k\mu}(\eq)$,
and the rate function $\rate_{k}$ is defined as in \eqref{eq:rate}.
} 

In particular, if $\theta_{k}'(0)$ is finite, convergence occurs in finite time.
\end{proposition}

\begin{Proof}
Pick some $\eps>0$ and let $U$ be a neighborhood of $\eq$ in $\strat$ such that $\payv_{k\alpha_{k}^{\ast}}(x) - \payv_{k\mu}(x)$ remains (multiplicatively) within $(1\pm\eps)$ of $\payv_{k\alpha_{k}^{\ast}}(\eq) - \payv_{k\mu}(\eq) = \delta_{k\mu}$ for all $x\in U$ and for all $\mu\in\act_{k}^{\ast} = \act_{k}\exclude\{k\alpha_{k}^{\ast}\}$.
If $x(0)$ is sufficiently close to $\eq$, $x(t)$ will remain in $U$ for all $t\geq0$ by Theorem \ref{thm:folk}.
\revised{%
Hence, the same reasoning as in the proof of Proposition \ref{prop:dom-rate} yields
\begin{equation}
x_{k\mu}(t)
	\leq \rate_{k}(c_{k\mu} - (1 - \eps) \tempinv_{k} \delta_{k\mu} t),
\end{equation}
and our assertion follows by summing over $\mu\in\act_{k}^{\ast}$.%
} 
\end{Proof}

\section{Time averages and the best response dynamics}
\label{sec:averages}

As is well known, the replicator dynamics \eqref{eq:RD} do not converge to equilibrium when the game's only Nash equilibrium is interior (for instance, as in \acl{MP} and generic zero-sum games).
On the other hand, if a replicator trajectory stays away from the boundary of the simplex in a $2$-player normal form game, its time-average converges to the Nash set of the game \textendash\ see e.g. \citeor[Chap.~7]{HS98}.

Under the reinforcement learning dynamics \eqref{eq:RL}, trajectories of play may enter and exit the boundary of the game's strategy space in perpetuity so the ``permanence'' criterion of staying a bounded distance away from $\bd(\strat)$ is no longer natural.
Instead, the conclusion about convergence of time averages can be reached by requiring that differences between the scores of each player's strategies remain bounded:

\begin{theorem}
\label{thm:average}
Let $x(t) = \choice(y(t))$ be an orbit of \eqref{eq:RL} in $\strat$ for a $2$-player game $\game$.
If the score differences $y_{k\alpha}(t) - y_{k\beta}(t)$ remain bounded for all $\alpha,\beta\in\act_{k}$, $k=1,2$, the time average $\widebar x(t) = t^{-1} \int_{0}^{t} x(s) \dd s$ of $x(t)$ converges to the set of Nash equilibria of $\game$.
\end{theorem}

As with the classic result on time averages, our proof relies on the linearity (as opposed to multilinearity) of each player's payoff function, a property specific to $2$-player games.

\begin{Proof}
By the definition of the dynamics \eqref{eq:RL} we have:
\begin{flalign}
y_{k\alpha}(t) - y_{k\beta}(t)
	&= c_{\alpha\beta}
	+ \int_{0}^{t} \left[ \payv_{k\alpha}(x(s)) - \payv_{k\beta}(x(s)) \right] \dd s
	\notag\\
	&= c_{\alpha\beta}
	+ t\cdot \left[ \payv_{k\alpha}(\widebar x(t) ) - \payv_{k\beta}(\widebar x(t)) \right],
\end{flalign}
where $c_{\alpha\beta} = y_{k\alpha}(0) - y_{k\beta}(0)$ and we used linearity to bring the integral inside the argument of $\payv_{k\alpha}$ and $\payv_{k\beta}$.
Thus, dividing by $t$ and taking the limit $t\to\infty$, we get
\begin{equation}
\lim_{t\to\infty}
	\left[ \payv_{k\alpha}(\widebar x(t)) - \payv_{k\beta}(\widebar x(t)) \right]
	= 0,
\end{equation}
where we have used the assumption that $y_{k\alpha} - y_{k\beta}$ is bounded.
Hence,
\revised{if $\eq$ is an $\omega$-limit point of $x(t)$},
we will have $\payv_{k\alpha}(\eq) = \payv_{k\beta}(\eq)$ for all $\alpha,\beta\in\act_{k}$, $k=1,2$, so $\eq$ must be a Nash equilibrium of $\game$.
Since $\strat$ is compact, the $\omega$-limit set of $x(t)$ is nonempty and our assertion follows.
\end{Proof}

\begin{remark}
To see how Theorem \ref{thm:average} implies the corresponding result for the replicator dynamics, simply note that $y_{k\alpha} - y_{k\beta} = \log x_{k\alpha} - \log x_{k\beta}$ under \eqref{eq:XL}.
Therefore, if $x(t)$ stays away from the boundary of the simplex, the requirement of Theorem \ref{thm:average} is fulfilled and we recover the standard result of \citeor{HS98}.
\end{remark}

The standard example of a $2$-player game that cycles under the replicator dynamics is the zero-sum game of \acf{MP} with payoff bimatrix:
\begin{equation}
U_{\textup{MP}} = 
\left(
\begin{array}{cc}
(1,-1)		& (-1,1)
\\
(-1,1)		& (1,-1)
\end{array}
\right).
\end{equation}
The game's unique minmax solution \citeorp{vN28} and unique Nash equilibrium is $\eq_{1} = \eq_{2} = (1/2,1/2)$,  and it is well known that the \acl{KL} divergence
\begin{equation}
\dkl(\eq,x)
	= \insum_{k} \insum_{\alpha}^{k} \eq_{k\alpha} \log \left(\eq_{k\alpha} \big/ x_{k\alpha}\right)
\end{equation}
is a constant of motion for \eqref{eq:RD} \citeorp{HS98}.
This implies that replicator trajectories always stay away from the boundary of $\strat$, so their time averages converge to the game's (unique) equilibrium.

Extending the above result to the reinforcement learning dynamics \eqref{eq:RL},
we have:

\begin{proposition}
\label{prop:zerosum}
Let $\game$ be a $2$-player zero-sum game \textup($\pay_{1} = -\pay_{2}$\textup)
\revised{%
that admits an interior equilibrium.
Then, the time average $\widebar x(t) = t^{-1} \int_{0}^{t} x(s) \dd s$ of every orbit $x(t) = \choice(y(t))$ of \eqref{eq:RL} in $\strat$ converges to the set of Nash equilibria of $\game$.
} 
\end{proposition}

\begin{Proof}
Let $\eq$ be an interior equilibrium of $\game$ and let $\fench_{h}(\eq,y) = \sum_{k=1,2} \big[ h_{k}(\eq_{k}) + h_{k}^{\ast}(y_{k}) - \langle{y_{k}}\vert{\eq_{k}}\rangle\big]$ denote the Fenchel coupling between $\eq$ and $y$ (cf. Appendix \ref{app:Bregman}).
Then, by Lemma \ref{lem:Fenchel-gradient}, we get:
\begin{flalign}
\label{eq:zerosum1}
\frac{d}{dt} \fench_{h}(\eq,y)
	&= \braket{\payv_{1}(x)}{x_{1} - \eq_{1}} + \braket{\payv_{2}(x)}{x_{2} - \eq_{2}}
	\notag\\
	&= \pay_{1}(x_{1},x_{2}) - \pay_{1}(\eq_{1},x_{2})
	+ \pay_{2}(x_{1},x_{2}) - \pay_{2}(x_{1},\eq_{2})
	= 0
\end{flalign}
on account of the game being zero-sum.
The above shows that $\fench_{h}(\eq,y(t))$ remains constant along \eqref{eq:RL}.
Proposition \ref{prop:Fenchel-bounded} then implies that $y_{\alpha}(t) - y_{\beta}(t)$ is bounded, so our assertion follows from Theorem \ref{thm:average}.
\end{Proof}

\begin{figure}[t]
\centering
\subfigure[Exponential reinforcement learning.]{
\label{subfig:averages.boundary}
\includegraphics[width=0.48\textwidth]{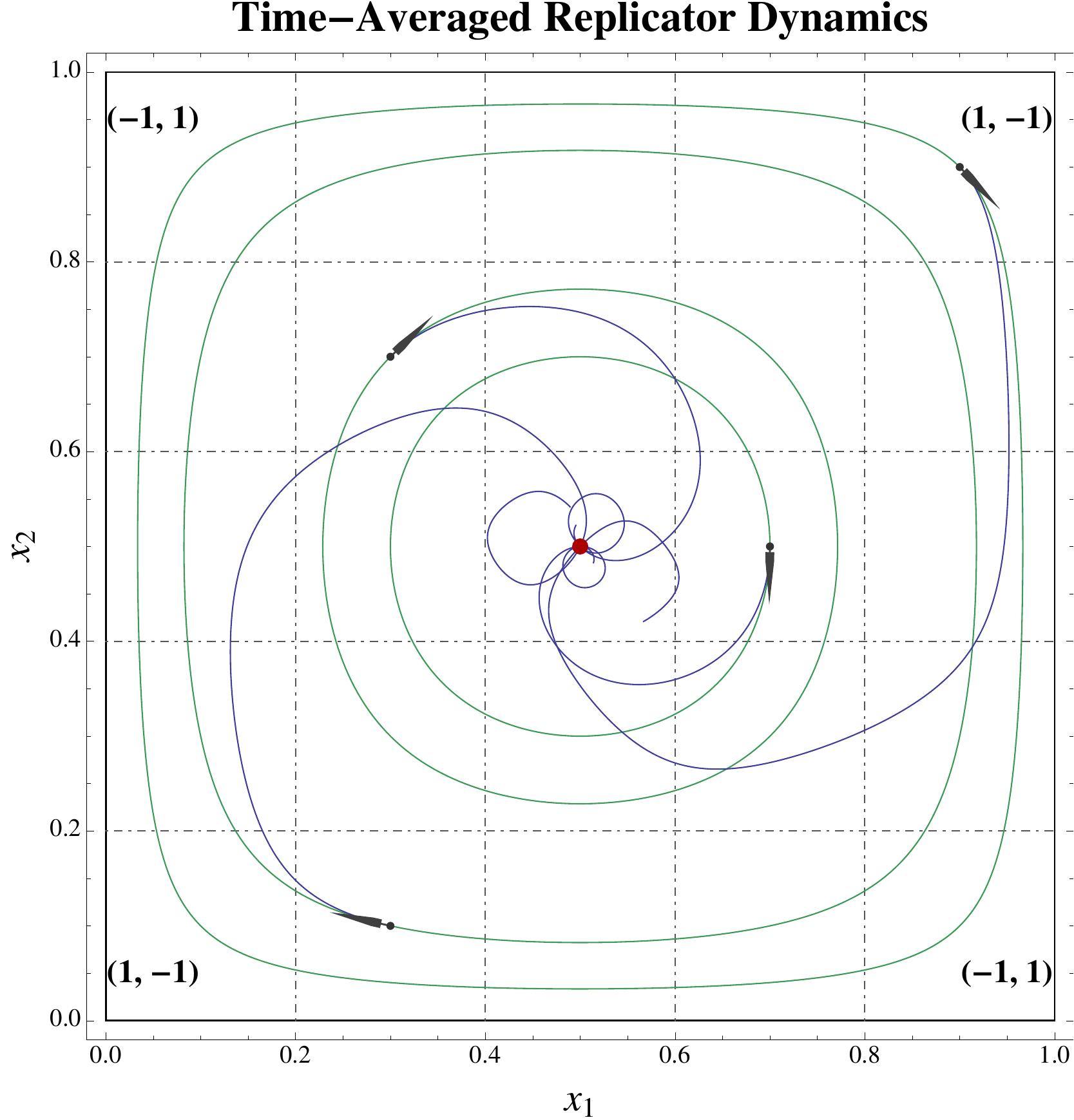}}
\hfill
\subfigure[Projected reinforcement learning.]{
\label{subfig:averages.long}
\includegraphics[width=0.48\textwidth]{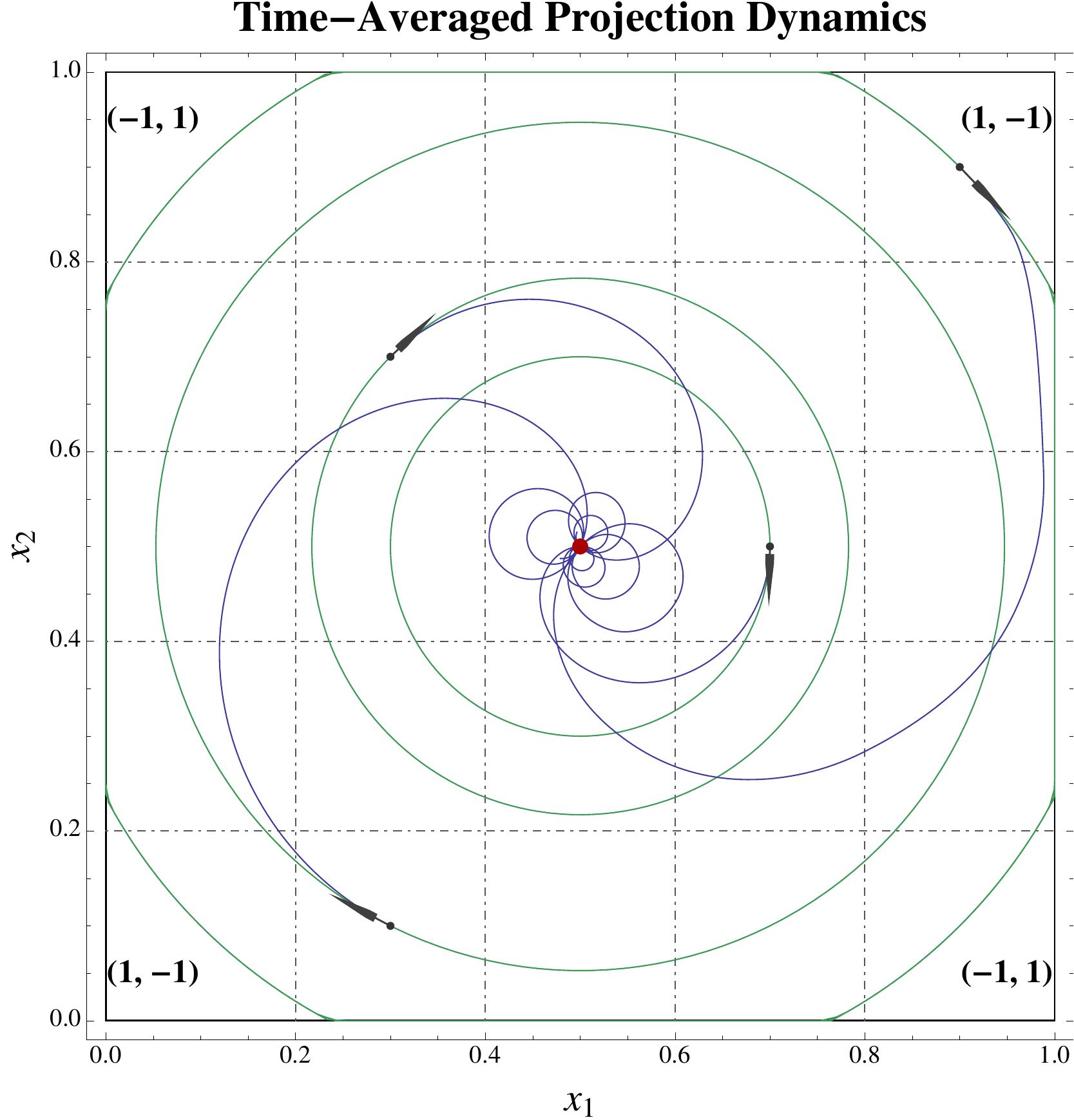}}
\caption{\small
Time averages in \acl{MP} under the exponential learning scheme \eqref{eq:XL} and the projected reinforcement learning dynamics \eqref{eq:PL}.
In both cases, solution trajectories (green) cycle, either avoiding the boundary $\bd(\strat)$ of $\strat$ (replicator) or attaining it infinitely often (projection);
on the other hand, their time averages (blue) converge to the game's equilibrium.
In both cases, the dynamics have been integrated over the same time horizon, showing that \eqref{eq:PL} evolves significantly faster than \eqref{eq:XL}.}
\label{fig:averages}
\end{figure}

In the case of the replicator dynamics, a heuristic explanation for the above result is that time averages of replicator trajectories in $2$-player games exhibit the same long-run behavior as the best response dynamics of \citeor{GM91}
\begin{equation}
\label{eq:BRD}
\tag{BRD}
\dot x_{k}
	\in \BREP_{k}(x) - x_{k}.
\end{equation}
with $\BREP_{k}(x) = \arg\max_{x_{k}'\in\strat_{k}} \langle{\payv_{k}(x)}\vert{x_{k}'}\rangle \subseteq\simplex$
denoting the \emph{best response correspondence} of player $k$.
\citeor{HSV09} showed that the $\omega$-limit set $\Omega$ of a time-averaged replicator orbit 
 is \emph{internally chain transitive} under \eqref{eq:BRD}: any two points $x,y\in\Omega$ may be joined by a piecewise continuous curve (a ``chain'') consisting of arbitrarily long pieces of orbits in $\Omega$ broken by arbitrarily small jump discontinuities (see \citeor{BHS05} for the precise definition).

As it turns out, this property extends verbatim to the learning scheme \eqref{eq:RL}:

\begin{theorem}
\label{thm:BRD}
Let $x(t) = \choice(y(t))$ be an orbit of \eqref{eq:RL} in $\strat$ for a $2$-player game $\game$.
Then, the $\omega$-limit set of the time average $\widebar x(t)$ of $x(t)$ is internally chain transitive under the best reply dynamics \eqref{eq:BRD}.
\end{theorem}

The proof of Theorem \ref{thm:BRD} follows closely that of \citeor[Proposition~5.1]{HSV09} and relies on the following proposition \textendash\ itself a generalization of (and proved in the same way as) Proposition 4.2 in \citeor{HSV09}:

\begin{proposition}
\label{prop:BRD-tracking}
Let $x(t) = \choice(y(t))$ be an orbit of \eqref{eq:RL} for a $2$-player game $\game$.
Then, $x_{k}(t)$ lies within $\delta(t)$ of $\BREP_{k}(\widebar x(t))$ in the uniform norm on $\strat$,
\revised{with $\delta(t)\to 0$ as $t\to\infty$}.
\end{proposition}

\begin{Proof}
From the definition of the dynamics \eqref{eq:RL}, and using the fact that $\payv_{k}$ is linear in 2-player games, we readily obtain: 
\begin{equation}
\score_{k}(t)
	= \score_{k}(0) + \int_{0}^{t} \payv_{k}(x(s)) \dd s
	= \score_{k}(0) + t \cdot \payv_{k}(\widebar x(t)).
\end{equation}
Hence, $x_{k}(t) = \choice_{k}(y_{k}(t))$ is the (unique) maximizer of the strictly concave problem:
\begin{equation}
\begin{aligned}
\text{maximize}
	\quad
	&\braket{\payv_{k}(\widebar x(t))}{x_{k}'} + t^{-1}\left[ \braket{y_{k}(0)}{x_{k}'} - h(x_{k}') \right],\\
\text{subject to}
	\quad
	&x_{k}'\in\strat_{k}.
\end{aligned}
\end{equation}
Therefore, with $\score(0)/t \to 0$ as $t\to\infty$ and $h$ finite and continuous on $\strat_{k}$, the maximum theorem of \citeor[p.~116]{Ber97} shows that $x_{k}(t)$ lies within a vanishing distance of $\BREP_{k}(\widebar x(t)) = \argmax_{x_{k}'\in\strat_{k}} \braket{\revised{\payv_{k}(\widebar x(t))}}{x_{k}'}$, as claimed.
\end{Proof}

\begin{Proof}[Proof of Theorem \ref{thm:BRD}]
Following \citeor{HSV09}, differentiate $\widebar x(t)$ to obtain
\begin{equation}
\frac{d}{dt} \widebar x(t)
	= \frac{x(t)}{t} - \frac{1}{t^{2}} \int_{0}^{t} x(s) \dd s
	=\frac{1}{t}\left(x(t) - \widebar x(t)\right).
\end{equation}
After changing time to $\tau = \log t$, this expression gives $\frac{d}{d\tau} \widebar x = x - \widebar x$, so Proposition \ref{prop:BRD-tracking} shows that $\widebar x(\tau)$ tracks a perturbed version of the best reply dynamics \eqref{eq:BRD} in the sense of \citeor[Definition II]{BHS05}.
Our assertion then follows from Theorem 3.6 in \citeor{BHS05}.
\end{Proof}

We close this section with some easy corollaries of Theorem \ref{thm:BRD}:
first, if the time average of a solution orbit $x(t) = \choice(y(t))$ of \eqref{eq:RL} converges,
\revised{its limit must be a Nash equilibrium};
second, if \eqref{eq:BRD} is globally attracted to some \revised{$\eq\in\strat$}, the time averages of \eqref{eq:RL} also converge to $\eq$.
These conclusions can be seen as generalizations of the corresponding statements in \citeor{HSV09} for \emph{interior} replicator orbits.
In our more general setting however, these conclusions hold for \emph{every} orbit of \eqref{eq:RL} in $\strat$, even those that enter and then leave $\bd(\strat)$  \textendash\ e.g.~as in the \acl{MP} example of Fig.~\ref{fig:averages}.

\section{Learning without a penalty function}
\label{sec:best}

We conclude this paper by considering the following question:
what happens if players use the \emph{exact} argmax correspondence $\brep_{k}(y_{k}) = \argmax_{x_{k}\in\strat_{k}} \braket{y_{k}}{x_{k}}$ of \eqref{eq:BR} as a choice map in \eqref{eq:RL}?

Given that $\brep_{k}\from \dual_{k} \rightrightarrows \strat_{k}$ is multi-valued, the resulting \ac{URL} process is defined via the inclusion:
\begin{equation}
\label{eq:BRL-int}
\tag{URL}
\begin{aligned}
\score_{k}(t)
	&= \score_{k}(0) + \int_{0}^{t} \payv_{k}(x(s)) \dd s,\\
x_{k}(t)
	&\in\brep_{k}(\score_{k}(t)).
\end{aligned}
\end{equation}
Expressing this in differential form yields the differential inclusion
\begin{equation*}
\dot\score_k(t)
	\in \payv_{k}(\revised{\brep(\score(t))}).
\end{equation*}

This unpenalized version of \eqref{eq:RL} is equivalent to the process of correlated \ac{CFP} (\citeor{FL98}),
in which players best respond to the time average of their opponents' joint past play.%
\footnote{We are grateful to an anonymous referee for pointing out this connection.}
To define this process formally, recall first that a \emph{correlated strategy} $\corr = (\corr_{\alpha_{1},\dotsc,\alpha_{N}}) \in\Delta(\act)$ is a distribution on pure strategy profiles $(\alpha_{1},\dotsc,\alpha_{N}) \in\act = \prod_{k}\act_{k}$.
Writing $\act_{-k} \equiv \prod_{\ell\neq k} \act_{\ell}$ for the set of action profiles of \revised{player $k$'s opponents},
let $\corr_{-k} \in \simplex(\act_{-k})$ denote the associated marginal distribution of $\corr$, and let $\payv_{k}^{c}(\corr)$ denote the payoff vector of player $k$ against this marginal distribution:
\begin{equation}
\label{eq:payvc}
\payv_{k\alpha_k}^{c}(\corr)
	= \insum_{\alpha_{-k}'\in\act_{-k}} \pay_{k}(\alpha_k,\alpha'_{-k}) \; \corr_{-k\alpha'_{-k}}.
\end{equation}
\revised{Finally}, given a strategy profile $x(t) = (x_{1}(t),\dotsc,x_{N}(t))$ at time $t$, let $\corr(t) = \bigotimes_{k} x_{k}(t)$ be its representation as a correlated strategy, i.e.
\begin{equation}\label{eq:MSPtoCorr}
\corr_{\alpha_{1},\dotsc, \alpha_{N}}(t)
	\equiv \inprod_{k} x_{k\alpha_k}(t).
\end{equation}

To define fictitious play in continuous time, we initialize the players' strategies $x_{k}(t)$ arbitrarily over some interval of time $[-\tau, 0)$, and thereafter set:
\begin{equation}
\label{eq:CFP}
\tag{CFP}
\begin{aligned}
\widebar \chi(t)
	&=\frac1{\tau+t} \int_{-\tau}^{t} \corr(s) \dd s,
	\\
x_{k}(t)
	&\in \brep_{k}(\payv^c_{k}(\widebar\corr(t))).
\end{aligned}
\end{equation}
Because the independent joint behavior at time $t$ is represented as a correlated strategy $\corr(t)$, the time-averaged joint behavior $\widebar \chi(t)$ typically exhibits correlation.%
\footnote{For instance, if the players first play the pure strategy profile $(\alpha_{1},\dotsc,\alpha_{N})$ and subsequently a pure strategy profile $(\alpha_{1}',\dotsc,\alpha_{N}')$ with $\alpha_k \ne \alpha_k'$ for all $k$, then $\widebar \chi(t)$ will be a joint distribution that puts all mass on these two strategy profiles.  This $\widebar \chi(t)$ cannot be represented as a mixed strategy profile, as it is not a product distribution.}

To show the equivalence of \eqref{eq:CFP} with the unpenalized reinforcement learning process \eqref{eq:BRL-int}, we use the linearity of $\payv_{k}^{c}$
\revised{and the identity $\payv^c_{k}(\corr(s))\equiv \payv_{k}(x(s))$}
to express $\payv_{k}^{c}(\widebar\corr(t))$ in terms of the score variables $\score_k(t)$
\revised{from \eqref{eq:BRL-int}}.
Defining the initial score as 
\begin{equation}\label{eq:CFPInit}
\score_k(0) = \int_{-\tau}^{0} \payv_{k}(x(s)),
\end{equation}
we obtain
\begin{flalign}
\payv^c_{k}(\widebar\corr(t))
	&= \payv_{k}^{c} \left( \frac1{\tau+t} \int_{-\tau}^{t} \corr(s) \dd s \right)
	\notag\\
	&= \frac1{\tau+t} \int_{-\tau}^{t} \payv^c_{k}(\corr(s)) \dd s
	\notag\\
	&=\frac1{\tau+t} \left( \score_k(0) +\int_{0}^{t} \payv_{k}(x(s)) \dd s \right)
	\notag\\
	&=\frac1{\tau+t} \,\score_k(t).
\end{flalign}
Since the argmax correspondence $\brep_k$ is scale-invariant, we can then write 
\begin{equation}
x_{k}(t)
	\in \brep_{k}(\payv_{k}^{c}(\widebar\corr(t)))
	= \brep_{k} \left( \frac{1}{\tau+t}\, \score_{k}(t) \right)
	= \brep_{k}(\score_{k}(t)),
\end{equation} 
thus recovering \eqref{eq:BRL-int}.%
\footnote{This last step would fail for our original process \eqref{eq:RL}, since the regularized argmax $\choice_k$ is not scale invariant.}

The only difference between \eqref{eq:CFP} and \eqref{eq:BRL-int} is in their allowable initial conditions:
Eq.~\eqref{eq:CFPInit} shows that \eqref{eq:CFP} can only generate initial score vectors $\score_{k}(0)$ corresponding to aggregate payoffs from an initial period of play.
%
%
Still, the definition of $\score_{k}(t)$ in \eqref{eq:BRL-int} implies that  as $t$ increases, the averaged score
\revised{$\frac1{t}\score_{k}(t)$ approaches the set of initial conditions of form \eqref{eq:CFPInit}.
In combination with the scale invariance of $\brep_k$, this constrains how different the asymptotic behaviors of \eqref{eq:BRL-int} and \eqref{eq:CFP} can be.}

In 2-player games, each player has just one opponent, so there is no need to account for correlation in the opponents' choices over time.  In this case we can express \eqref{eq:CFP} as
\begin{equation}
\begin{aligned}
{\widebar x}_{k}(t)
	&=\frac1{\tau+t} \int_{-\tau}^{t} x_k(s) \dd s,
	\\
x_{k}(t)
	&\in \brep_{k}(\payv_{k}(\widebar x(t))).
\end{aligned}
\end{equation}
Differentiating then yields
\begin{equation}
\frac{d}{dt}\widebar x_{k}(t)
	\in \frac1{\tau+t} \left(\brep_{k}(\payv_{k}(\widebar x(t))) - \widebar x(t) \right)\\ 
	= \frac1{\tau+t} \left(\BREP(\widebar x(t)) - \widebar x(t) \right).
\end{equation}
Thus, up to a time change, the evolution of each player's time-averaged play under two-player \eqref{eq:BRL-int} and \eqref{eq:CFP} follows the best response dynamics \eqref{eq:BRD}.
Consequently, standard results on these dynamics (see e.g.~\citeor{Hof95}) imply that analogues of certain properties of \eqref{eq:RL} still hold for \eqref{eq:BRL-int} when initial scores are of the form \eqref{eq:CFPInit}: 
\begin{enumerate}
\item
Dominated strategies are never chosen. 
\item
There exists a stationary trajectory $x(t)$ of \eqref{eq:BRL-int} such that $x(t) = \eq$ (or $\widebar x(t) = \eq$) for all $t\geq0$, if and only if $\eq$ is a Nash equilibrium.
\item
If $\eq$ is a strict equilibrium, then it attracts an open set of nearby initial conditions $\widebar x(0)$.
Alternatively, there is an open set of initial conditions $y(0)$ with $\brep(y(0))$ close enough to $\eq$ such that $\widebar x(t)$ and $x(t)$ converge to $\eq$.
\item
In $2$-player zero-sum games, $\widebar x(t)$ converges to the game's set of Nash equilibria.
\end{enumerate}
It is easy to establish versions of these results that allow for arbitrary initial score vectors.
The only amendment needed is to claim (1), to say that dominated strategies cease to be chosen after some finite time interval.

With three or more players, \eqref{eq:BRL-int} and \eqref{eq:CFP} are no longer equivalent to \eqref{eq:BRD}:
as we have seen, the former two processes incorporate the correlation which comes from averaging over time while the latter process does not. 
In this setting, \eqref{eq:CFP} does not define a convex-valued differential inclusion, so the basic properties of its solutions do not follow from standard results.%
\footnote{Whenever $\widebar\corr(t)$ is a correlated strategy against which multiple players have multiple best responses, the set of allowable choices of $\corr(t)$ in \eqref{eq:CFP} is a set of correlated strategies corresponding to independent randomizations by these players.
This is not a convex subset of $\Delta(\act)$ (cf.~\citeor{VZ13}).
\revised{For similar reasons, if each player's score vector $\score_k(t)$ admits multiple best responses, the set of feasible vectors of motion $\dot\score(t)$ under \eqref{eq:BRL-int} need not be convex.}
For basic properties of differential inclusions, see \citeor{AC84}.}
Still, properties (1)\textendash(3) above are easily established for the process \eqref{eq:CFP}, and therefore hold for \eqref{eq:BRL-int} as well.

\begin{figure}[t]
\centering
\subfigure{
\label{fig:best-logit}
\includegraphics[width=0.48\textwidth]{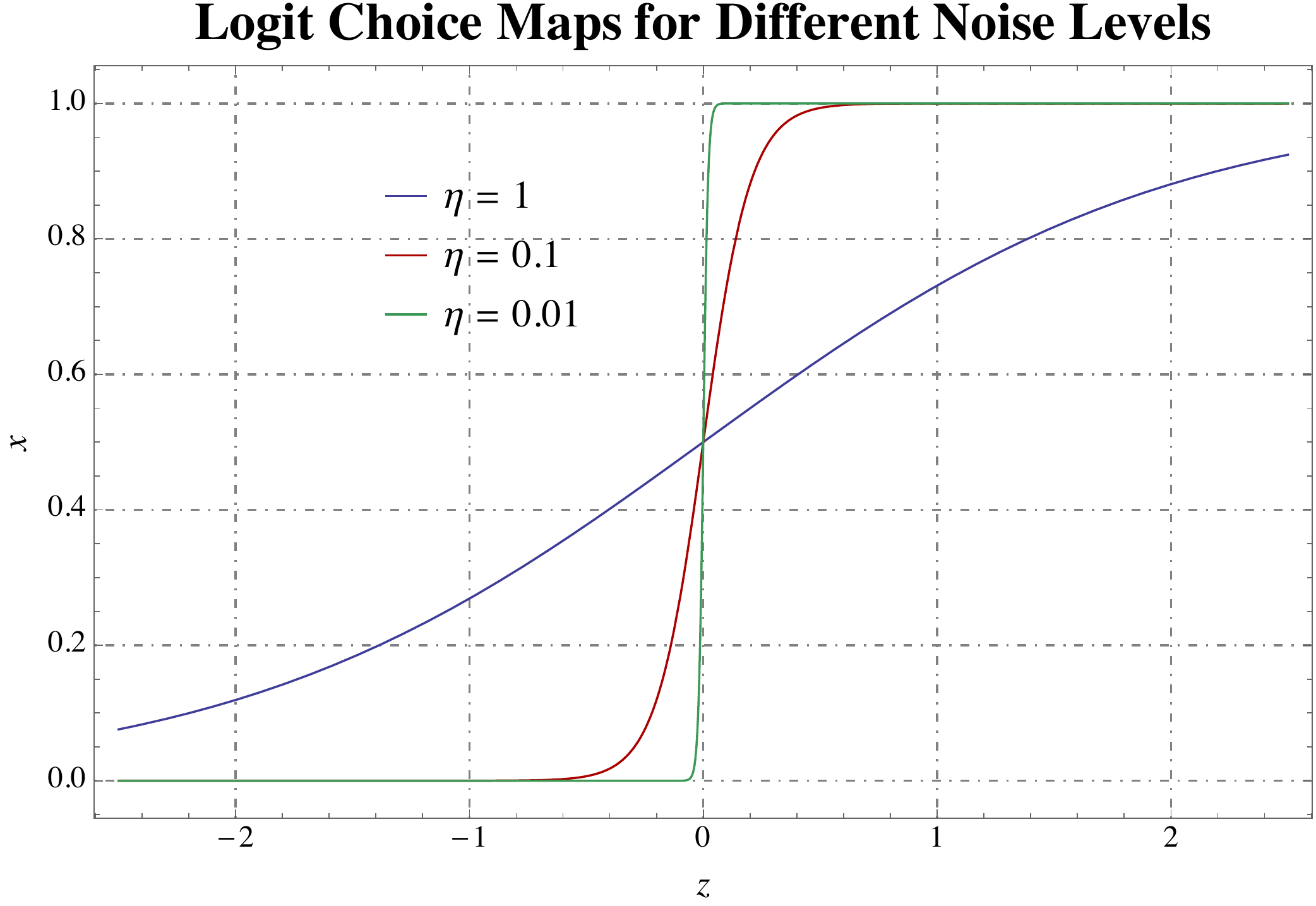}}
\hfill
\subfigure{
\label{fig:best-proj}
\includegraphics[width=0.48\textwidth]{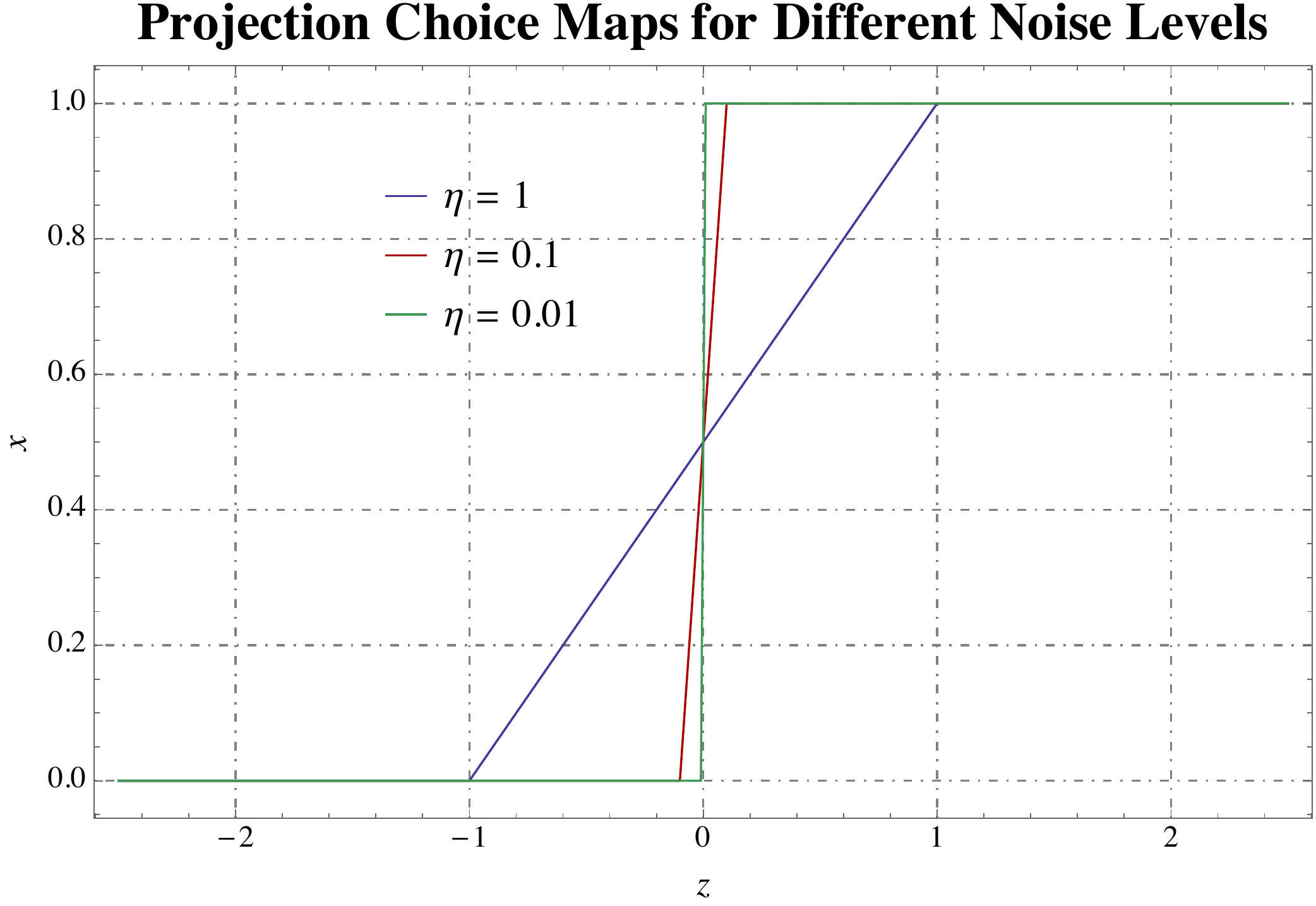}}
\caption{\small
Logit and projection choice maps for different noise levels.
We consider a player with two actions $\act=\{0,1\}$ and $z = y_{1} - y_{0}$ denotes the score difference between the two.
The player's choice probability $\choice_{1}(y)$ for different noise levels $\temp = \tempinv^{-1}$ is given by \eqref{eq:choice-logit2} for the logit case, and by $\min\{1,\max\{0, (1 + \temp^{-1} z)/2\}\}$
for the projection case.}
\label{fig:best}
\end{figure}


The rationality properties of \eqref{eq:BRL-int} can also be explained by considering the reinforcement learning scheme \eqref{eq:RL} with a very small noise level $\temp \approx 0$ as in \eqref{eq:RL-temp}.
In Fig.~\ref{fig:best}, we focus on the entropic and quadratic penalties \eqref{eq:penalty-Gibbs} and \eqref{eq:penalty-Eucl}, and we plot the induced regularized best response maps for different noise levels $\temp$ when the player has two pure strategies.
Up to a time shift, these plots also describe the behavior of the solution orbits of \eqref{eq:RL-temp} for a single-player game 
with fixed payoff difference $\payv_{1}-\payv_{0} \equiv 1$ (simply note that the player's score difference $z$ satisfies $\dot z = \payv_{1}-\payv_{0} = 1$).
In the logit case, the weight on the dominated strategy vanishes at an exponential rate,  and this rate of decay increases as the noise level $\temp$ tends to $0$ (cf.~Table \ref{tab:rates});
nevertheless, the dominated strategy is played with positive probability for all time.
In the projection case however, the player assigns positive weight on both actions only within a window whose width goes to $0$ as the noise level $\temp$ tends to $0$;
thus, if the player starts with a high score difference in favor of the dominated strategy, this difference will decrease linearly until it reaches that window, and the player will then transit sharply to the dominant strategy.
Either way, the behavior of the dynamics \eqref{eq:BRL-int} is recovered in the limit $\temp \to 0$.

\appendix

\section{Basic properties of choice maps}
\label{app:choice}

\numberwithin{remark}{section}

In this appendix, we prove some intuitive properties of choice maps that were invoked in the main text.
For simplicity, all results will be stated for the unit simplex $\simplex \equiv \simplex(\act)$ of $\primal=\R^{n}$ that is spanned by the index set $\act = \{1,\dotsc,n\}$.
We have:

\begin{proposition}
\label{prop:choice}
Let $h$ be a penalty function on $\simplex$ and let $\choice\from \dual \to \simplex$ be its induced choice map.
Then:
\begin{enumerate}
[i.]
\addtolength{\itemsep}{1pt}

\item
$\choice(y') = \choice(y)$ for all $y,y'\in \dual$ such that $y' - y \propto (1,\dotsc,1)$.

\item
$\choice(y-t\bvec_{\beta}) = \choice(y)$ for all $t\geq0$ and for all $y\in \dual$ such that $\choice_{\beta}(y) = 0$.

\item
$\choice_{\alpha}(y)\to0$ whenever $y_{\alpha} - y_{\beta} \to -\infty$ for some $\beta\neq\alpha$.
\end{enumerate}
\end{proposition}

The first part of Proposition \ref{prop:choice} shows that modifying all payoffs by the same amount does not change relative advantages between strategies, so choice probabilities remain invariant along $(1,\dotsc,1)$.
The second part shows that choice probabilities also remain unchanged if one reduces the payoff of an action that already has a strategy share of $0$.
Finally, the last part states that the strategy share of an action becomes vanishingly small when the payoff of said action is at a great relative disadvantage to that of another action.

\begin{Proof}[Proof of Proposition \ref{prop:choice}]
For the first part, assume that $y_{\alpha}' = y_{\alpha} + c$ for some $c\in\R$ and for all $\alpha\in\act$.
Then:
\begin{flalign}
\choice(y')
	&\txs
	= \argmax_{x\in\simplex} \big\{\braket{y}{x} + \insum_{\beta} x_{\beta}\cdot c - h(x)\big\}
	\notag\\
	&\txs
	= \argmax_{x\in\simplex} \{\braket{y}{x} - h(x)\}
	= \choice(y).
\end{flalign}

For the second part, the first order \acf{KKT} conditions for the regularized problem \eqref{eq:softmax} give
\begin{equation}
\label{eq:KKT}
y_{\beta} - \left.\frac{\pd h}{\pd x_{\beta}}\right\vert_{\choice(y)}
	= \mu - \nu_{\beta}
	\qquad
	(\beta = 1,\dotsc,n),
\end{equation}
where $\mu\in\R$ is the Lagrange multiplier of the equality constraint $\insum_{\beta} x_{\beta} = 1$ and $\nu_{\alpha} \geq 0$ is the complementary slackness multiplier of the inequality constraint $x_{\alpha}\geq0$ \textendash\ i.e. $\nu_{\alpha} = 0$ for all $\alpha\in\supp(\choice(y))$.
Hence, if we let $y_{\alpha}' = y_{\alpha} - t\delta_{\alpha\beta}$ for some $\beta\notin\supp(\choice(y))$ and set $\mu' = \mu$ and $\nu_{\alpha}' = \nu_{\alpha} + t \delta_{\alpha\beta} \geq 0$,
$\choice(y)$ also satisfies the \ac{KKT} conditions of \eqref{eq:softmax} for $y'$ with Lagrange multipliers $\mu'$ and $\nu_{\alpha}'$.
Since $h$ is strictly convex, it follows that $\choice(y') = \choice(y)$, as claimed.

As for the last part, let $y_{j}$ be a sequence in $\dual$ such that $y_{j,\alpha} - y_{j,\beta} \to -\infty$, set $x_{j} = \choice(y_{j})$, and assume (by descending to a subsequence if necessary) that $x_{j,\alpha} > \eps >0$ for all $n$.
By definition, we then have
\begin{flalign}
\label{eq:Qcomp1}
\braket{y_{j}}{x_{j}} - h(x_{j})
	\geq \braket{y_{j}}{x'} - h(x')
\end{flalign}
for all $x'\in\simplex$.
Therefore, if we set $x_{j}' = x_{j} + \eps(\bvec_{\beta} - \bvec_{\alpha})$, we readily get
\begin{equation}
\label{eq:Qcomp2}
\eps (y_{j,\alpha} - y_{j,\beta})
	\geq h(x_{j}) - h(x_{j}')
	\geq - (h_{\max} - h_{\min}),
\end{equation}
which contradicts our original assumption that $y_{j,\alpha} - y_{j,\beta} \to -\infty$.
With $\simplex$ compact, the above shows that $x_{\alpha}^{\ast} = 0$ for any limit point $x^{\ast}$ of $x_{j}$, i.e. $\choice_{\alpha}(y_{j})\to0$.
\end{Proof}

\section{Calculations for the R\'enyi dynamics}
\label{app:Renyi}

Here we provide the calculations leading to the Rényi dynamics \eqref{eq:ReD}.
To begin with, if we drop the player index $k$, a straightforward differentiation of the Rényi penalty function \eqref{eq:penalty-Renyi} yields:
\begin{equation}
\label{eq:Renyi1}
\frac{\pd h}{\pd x_{\alpha}}
	= \frac{q}{q-1} \frac{x_{\alpha}^{q-1}}{\insum_{\gamma} x_{\gamma}^{q}}
	= \frac{\xi_{\alpha}}{q-1},
\end{equation}
and hence:
\begin{equation}
\label{eq:Renyi2}
\frac{\pd^{2} h}{\pd x_{\alpha} \pd x_{\beta}}
	= q \delta_{\alpha\beta} \frac{x_{\alpha}^{q-2}}{\insum_{\gamma} x_{\gamma}^{q}}
	- \frac{q^{2}}{q-1}
	\frac{x_{\alpha}^{q-1}}{\insum_{\gamma} x_{\gamma}^{q}}
	\frac{x_{\beta}^{q-1}}{\insum_{\gamma} x_{\gamma}^{q}}
	= \delta_{\alpha\beta} \frac{\xi_{\alpha}}{x_{\alpha}}
	+ \frac{1}{1-q} \xi_{\alpha} \xi_{\beta},
\end{equation}
where, for simplicity, we have set:
\begin{equation}
\label{eq:xi}
\xi_{\alpha}
	= \frac{q x_{\alpha}^{q-1}}{\insum_{\gamma} x_{\gamma}^{q}}.
\end{equation}
Thus, letting $g_{\alpha\beta} = \pd_{\alpha} \pd_{\beta} h$, the next step is to calculate the inverse Hessian matrix $g^{\alpha\beta} = \hess(h(x))_{\alpha\beta}^{-1}$.
To that end, we claim that:
\begin{equation}
\label{eq:Renyi-inverse}
g^{\alpha\beta}
	= -x_{\alpha} x_{\beta} + \delta_{\alpha\beta} \frac{x_{\alpha}}{\xi_{\alpha}}.
\end{equation}
Indeed, by a straightforward verification, we get:
\begin{flalign}
\insum_{\beta} g^{\alpha\beta} g_{\beta\gamma}
	&= \insum_{\beta}
	\left( -x_{\alpha} x_{\beta} + \delta_{\alpha\beta} \frac{x_{\alpha}}{\xi_{\alpha}} \right)
	\left( \frac{1}{1-q} \xi_{\beta} \xi_{\gamma} + \delta_{\beta\gamma} \frac{\xi_{\beta}}{x_{\beta}} \right)
	\notag\\
	&= -\frac{1}{1-q} x_{\alpha} \xi_{\gamma} \insum_{\beta} x_{\beta} \xi_{\beta}
	- x_{\alpha} \xi_{\gamma}
	+ \frac{1}{1-q} x_{\alpha} \xi_{\gamma}
	+ \delta_{\alpha\gamma} \frac{\xi_{\alpha}}{x_{\alpha}} \frac{x_{\gamma}}{\xi_{\gamma}}
	\notag\\	
	&= \delta_{\alpha\gamma},
\end{flalign}
where we have used the fact that $\insum_{\beta} x_{\beta} \xi_{\beta} = q$ (by definition).

On that account, we obtain:
\begin{equation}
\label{eq:Renyi-normal}
g^{\alpha}
	= \insum_{\beta} g^{\alpha\beta}
	= -x_{\alpha} + \frac{x_{\alpha}}{\xi_{\alpha}}
\end{equation}
and
\begin{equation}
\label{eq:Renyi-norm}
G
	= \insum_{\alpha} g^{\alpha}
	= -1 + \insum_{\alpha} \frac{x_{\alpha}}{\xi_{\alpha}}
	= -1 + \frac{1}{q} \insum_{\alpha} x_{\alpha}^{q} \cdot \insum_{\alpha} x_{\alpha}^{2-q}.
\end{equation}
Thus, letting $S_{q} \equiv q^{-1} \insum_{\alpha} x_{\alpha}^{q} \cdot \insum_{\alpha} x_{\alpha}^{2-q} = 1 + G$, some more algebra yields:
\begin{flalign}
\label{eq:Renyi3}
g^{\alpha\beta} - \frac{g^{\alpha} g^{\beta}}{G}
	&= -x_{\alpha} x_{\beta} + \delta_{\alpha\beta} \frac{x_{\alpha}}{\xi_{\alpha}}
	+ \frac{x_{\alpha} x_{\beta} - x_{\alpha} x_{\beta} \xi_{\beta}^{-1} - x_{\alpha} x_{\beta} \xi_{\alpha}^{-1} + x_{\alpha}x_{\beta} \xi_{\alpha}^{-1}\xi_{\beta}^{-1}}{1 - S_{q}}
	\notag\\
	&= \frac{x_{\alpha}}{\xi_{\alpha}} \delta_{\alpha\beta}
	+ x_{\alpha} \frac{S_{q} - \xi_{\alpha}^{-1}}{1 - S_{q}} x_{\beta}
	- x_{\alpha} \frac{1 - \xi_{\alpha}^{-1}}{1- S_{q}} \frac{x_{\beta}}{\xi_{\beta}},
\end{flalign}
and \eqref{eq:ReD} follows immediately from Proposition \ref{prop:RLD}.

We are left to show that the limit of \eqref{eq:ReD} as $q\to1$ is $x_{\alpha}\big( \payv_{\alpha} - \insum_{\beta} x_{\beta} \payv_{\beta} \big)$, i.e., that it boils down to the RHS of the replicator dynamics \eqref{eq:RD}.
To do so, given that $\xi_{\alpha}\to1$ as $q\to1$, it suffices to show that
\begin{equation}
\label{eq:Renyi-limit1}
\lim_{q\to1} \left[
	\frac{S_{q} - \xi_{\alpha}^{-1}}{1 - S_{q}}
	- \frac{1}{\xi_{\beta}} \frac{1 - \xi_{\alpha}^{-1}}{1- S_{q}}
	\right]
	= -1,
\end{equation}
However, after discarding terms that tend to $1$ as $q\to1$, the above limit may be written as:
\begin{equation}
\label{eq:Renyi-limit2}
\lim_{q\to 1}
	\frac{\insum_{\gamma} x_{\gamma}^{2-q} - x_{\alpha}^{1-q} - x_{\beta}^{1-q} + q^{-1} x_{\alpha}^{1-q} x_{\beta}^{1-q} \insum_{\gamma} x_{\gamma}^{q}}
	{q - \insum_{\gamma} x_{\gamma}^{q} \cdot \insum_{\gamma} x_{\gamma}^{2-q}}.
\end{equation}
Hence, by using de l'Hôpital's rule, the limit \eqref{eq:Renyi-limit2} is equal to:
\begin{equation}
\frac
	{-\insum_{\gamma} x_{\gamma} \log x_{\gamma} + \log x_{\alpha} + \log x_{\beta} - 1 - \log x_{\alpha} - \log x_{\beta} + \insum_{\gamma} x_{\gamma} \log x_{\gamma}}
	{1 - \insum_{\gamma} x_{\gamma} \log x_{\gamma} + \insum_{\gamma} x_{\gamma} \log x_{\gamma}}
	= -1,
\end{equation}
as claimed.

\section{Bregman divergences and the Fenchel coupling}
\label{app:Bregman}

\numberwithin{remark}{section}

In this appendix, we introduce Bregman divergences and the Fenchel coupling, and we discuss their basic properties.

As before, let $h\from\simplex\to\R$ be a penalty function on the unit simplex $\simplex$ of $\primal \equiv \R^{n}$.
For convenience, we will treat $h$ as an extended-real-valued function $h\from \primal\to\R\cup\{+\infty\}$ by setting $h(x) = +\infty$ for all $x\in \primal\exclude\simplex$.
The \emph{subdifferential} of $h$ at $x\in \primal$ is then defined as $\subd h(x) = \{y \in \dual: h(x') \geq h(x) + \braket{y}{x' - x}\;\text{for all $x'\in \primal$}\}$ and we will say that $h$ is \emph{subdifferentiable at $x$} whenever $\subd h(x)$ is nonempty.
This is always the case if $x\in\intsimplex$, so we have $\intsimplex \subseteq \dom \subd h \equiv \{x\in\simplex: \subd h(x) \neq \varnothing\} \subseteq \simplex$.

A key tool in our analysis is the \emph{convex conjugate} $h^{\ast}\from \dual\to\R$ of $h$ defined as
\begin{equation}
\label{eq:conjugate}
h^{\ast}(y)
	= \max\nolimits_{x\in\simplex} \{\braket{y}{x} - h(x)\}.
\end{equation}
As it turns out, the choice map $\choice\from \dual \to\simplex$ induced by $h$ is simply the differential of $h^{\ast}$:

\begin{proposition}
\label{prop:choice-dh}
Let $h$ be a penalty function on $\simplex$.
The induced choice map $\choice\from \dual\to\simplex$ is Lipschitz and $\choice(y) = dh^{\ast}(y)$ for all $y\in \dual$.
\end{proposition}

\begin{Proof}
By Theorem 23.5 in \citeor{Roc70}, we readily obtain
\begin{equation}
\label{eq:dh}
x \in \subd h^{\ast}(y)
	\iff y \in \subd h(x)
	\iff x \in \argmax\nolimits_{x'\in\simplex} \{\braket{y}{x'} - h(x')\}.
\end{equation}
Since the last set only contains $\choice(y)$, we immediately obtain $\choice(y) = dh^{\ast}(y)$.
The Lipschitz property for $\choice$ then follows from the strong convexity of $h$ \textendash\ see e.g. \citeor{Nes09}.
\end{Proof}

\begin{remark}
\label{rem:imQ}
Proposition \ref{prop:choice-dh} is folklore in convex optimization \textendash\ see e.g. \citeor{HS02}, \citeor{Nes09}, \citeor{SS11}, \citeor{KM14} and many others.
Equation \eqref{eq:dh} also shows that the image of $\choice$ is precisely $\dom\subd h$, a fact which we use freely in the rest of this appendix.
\end{remark}

Given a basepoint $p\in\simplex$, the one-sided derivative
\begin{equation}
\label{eq:hder}
h'(x;p-x)
	= \lim_{t\to0^{+}} t^{-1} \left[h(x + t(p-x)) - h(x)\right]
\end{equation}
exists for all $x\in\simplex$ and is finite whenever $x$ lies in the relative interior of a face of $\simplex$ that also contains $p$.
With this in mind, we define the \emph{Bregman divergence} of $h$ as
\begin{equation}
\label{eq:Bregman}
\breg_{h}(p, x)
	= h(p) - h(x) - h'(x;p-x)
	\quad
	\text{for all $p,x\in\simplex$,}
\end{equation}
with $\breg_{h}(p,x)$ possibly attaining the value $+\infty$ if $h'(x;p-x) = -\infty$.%
\footnote{Usually, Bregman divergences are defined for $x\in\dom\subd h$ \textendash\ \citeor{Kiw97b} uses the notation $\breg_{h}'$ to distinguish \eqref{eq:Bregman} from the original definition of \citeor{Bre67}.
The ``\emph{raison d'\,être}'' of the more general definition \eqref{eq:Bregman} is that we often need to work with boundary points $x\in\bd(\simplex)$ with $\subd h(x) = \varnothing$.}
We then have:

\begin{proposition}
\label{prop:Bregman}
Let $h$ be a $K$-strongly convex penalty function on $\simplex$ and let $\simplex_{p}$ be the union of the relative interiors of the faces of $\simplex$ that contain $p$, i.e.
\begin{equation}
\label{eq:domp}
\simplex_{p}
	= \{x\in\simplex: \supp(x) \supseteq \supp(p)\}
	=\{x\in\simplex: \text{$x_{\alpha}>0$ whenever $p_{\alpha}>0$}\}.
\end{equation}
Then:
\begin{enumerate}
[i.]
\item
$\breg_{h}(p,x) < +\infty$ whenever $x\in\simplex_{p}$.

\item
$\breg_{h}(p,x)\geq0$ for all $x\in\simplex$ and $\breg_{h}(p,x) = 0$ if and only if $p=x$;
in particular:
\begin{equation}
\label{eq:Dbound}
\breg_{h}(p,x)
	\geq \frac{1}{2} K \norm{x - p}^{2}
	\quad
	\text{for all $x\in\simplex$.}
\end{equation}

\item
$\breg_{h}(p,x_{j}) \to \breg_{h}(p,x)$ whenever $x_{j}\to x$ in $\simplex_{p}$.
\end{enumerate}
\end{proposition}


\begin{Proof}
Let $z = p-x$.
If $x\in\simplex_{p}$, $h(x + tz)$ is finite and smooth for all $t$ in a neighborhood of $0$ so $h'(x;p-x)$ \textendash\ and hence $\breg_{h}(p,x)$ \textendash\ must be finite as well.

For Part (ii), positive-definiteness is a trivial consequence of strict convexity;
on the other hand, \emph{strong} convexity yields:
\begin{equation}
\label{eq:Dbound1}
h(x+tz)
	\leq t h(p) + (1-t) h(x) - \frac{1}{2} K t(1-t) \norm{z}^{2}.
\end{equation}
Moreover, with $h(x)$ finite, we also get $h(x+tz) \geq h(x) + t h'(x;z)$, so \eqref{eq:Dbound1} gives:
\begin{equation}
\label{eq:Dbound3}
h(x) + th'(x;z)
	\leq t h(p) +(1-t) h(x) - \frac{1}{2} K t(1-t) \norm{z}^{2}.
\end{equation}
After rearranging and dividing by $t$, the above becomes
\begin{equation}
\label{eq:Dbound4}
h(p) - h(x) - h'(x;z)
	\geq \frac{1}{2} K(1-t) \norm{z}^{2},
\end{equation}
so \eqref{eq:Dbound} is obtained by letting $t\to 0$.

Finally, if $x_{j}\to p$, Part (iii) follows from \citeor[Lemma~8.2]{Kiw97b}.
Otherwise, if $\lim_{j} x_{j} \neq p$, let $z_{j} = p - x_{j}$ and take $\eps>0$ such that $x_{j} + tz_{j}\in\simplex$ for all $t\in(-\eps,\eps)$ and for all sufficiently large $n$ (recall that $\simplex_{p}$ is relatively open in $\simplex$ so $x_{j}\in\simplex_{p}$ for large enough $n$).
Furthermore, let $f_{j}(t) = h(x_{j}+t z_{j})$ and let $f(t) = h(x+tz)$ for $t\in(-\eps,\eps)$;
since $f_{j}$ and $f$ are smooth, convex and $f_{j}\to f$ pointwise, we obtain $h'(x_{j};p-x_{j}) = f_{j}'(0)\to f'(0) = h'(x;p-x)$ by Theorem 25.7 in \citeor{Roc70}.
This concludes our proof.
\end{Proof}

Dually to the above, $h$ also induces a convex coupling $\fench_{h}\from \simplex\times \dual \to \R$ with
\begin{equation}
\label{eq:Bregman-dual}
\fench_{h}(p,y)
	= h(p) + h^{\ast}(y) - \braket{y}{p}
	\quad
	\text{for all $p\in\simplex$, $y\in \dual$.}
\end{equation}
This primal-dual coupling is (strictly) convex in both arguments and positive-semidefinite by Fenchel's inequality, so we call it the \emph{Fenchel coupling} between $p$ and $y$.
The next proposition establishes the duality relation between $\breg_{h}$ and $\fench_{h}$:

\begin{proposition}
\label{prop:Bregman-dual}
Let $h$ be a $K$-strongly convex penalty function on $\simplex$ and let $p\in\simplex$.
Then, $\fench_{h}(p,y) \geq \frac{1}{2} K \norm{\choice(y)-p}^{2}$ for all $y\in \dual$ and $\fench_{h}(p,y) \to 0$ if and only if $\choice(y) \to p$;
furthermore:
\begin{equation}
\label{eq:PD-duality}
\fench_{h}(p,y)
	= \breg_{h}(p,x)
	\quad
	\text{whenever $\choice(y) = x$ and $x\in\simplex_{p}$.}
\end{equation}
\end{proposition}

\begin{remark}
The first part of Proposition \ref{prop:Bregman-dual} is not implied by the second because $\im\choice=\dom\subd h$ is not necessarily contained in $\simplex_{p}$.
\end{remark}

\begin{Proof}
For the first part of the proposition, let $x = \choice(y)$ so that $h^{\ast}(y) = \braket{y}{x} - h(x)$.
Then:
\begin{equation}
\label{eq:Bregman-dual2}
\fench_{h}(p,y)
	= h(p) - h(x) - \braket{y}{p-x}.
\end{equation}
With $y\in\subd h(x)$, we also get
\begin{equation}
\label{eq:Pbound1}
h(x + t(p-x))
	\geq h(x) + t\braket{y}{p-x}
	\quad
	\text{for all $t\in[0,1]$},
\end{equation}
so, by combining \eqref{eq:Pbound1} with \eqref{eq:Dbound1} as in the proof of Prop.~\ref{prop:Bregman}, we obtain
\begin{equation}
\label{eq:Pbound2}
h(p) - h(x) - \braket{y}{p-x}
	\geq \frac{1}{2} K \norm{x-p}^{2},
\end{equation}
and our claim follows.

As for \eqref{eq:PD-duality}, if $x\in\simplex_{p}$, the function $h(x + t(p-x))$ is finite and smooth for all $t$ in a neighborhood of $0$ (simply note that $x+t(p-x)$ lies in the relative interior of a face of $\simplex$ for small $t$).
Thus, since $y\in\subd h(x)$ and $h$ admits a two-sided derivative at $x$ along $p-x$, we also have $h'(x;p-x) = \braket{y}{p-x}$, so \eqref{eq:PD-duality} follows from \eqref{eq:Bregman-dual2}.
\end{Proof}

Intuitively, Proposition \ref{prop:Bregman-dual} shows that the Fenchel coupling $\fench_{h}(p,y)$ between $p$ and $y$ measures the proximity of $\choice(y)$ to $p$;
as such, if $\fench_{h}(p,y)$ grows large, $\choice(y)$ must be moving away from $p$.
We formalize this as follows:

\begin{proposition}
\label{prop:divergence}
If $\fench_{h}(p,y_{j})\to+\infty$ for some sequence $y_{j}\in \dual$, the sequence $x_{j} = \choice(y_{j})$ has no limit points in $\simplex_{p}$;
in particular, $\liminf_{j\to\infty} \{x_{j,\alpha}: \alpha\in\supp(p)\} = 0$.
\end{proposition}

\begin{Proof}
By descending to a subsequence of $x_{j} = \choice(y_{j})$ if necessary, assume that $\lim_{j} x_{j} = \eq\in\simplex_{p}$.
Since $\simplex_{p}$ is relatively open in $\simplex$, we must eventually have $x_{j}\in\simplex_{p}$, so Proposition \ref{prop:Bregman} gives $\breg_{h}(p,x_{j})\to \breg_{h}(p,\eq) < +\infty$.
However, with $x_{j}\in\simplex_{p}$, Proposition \ref{prop:Bregman-dual} also yields $\breg_{h}(p,x_{j}) = \fench_{h}(p,y_{j})\to+\infty$, a contradiction.
\end{Proof}

The following result may be seen as a weak partial converse to the above:

\begin{proposition}
\label{prop:Fenchel-bounded}
Let $p\in\intsimplex$.
If $\fench_{h}(p,y_{j})$ is bounded, then $y_{j,\alpha} - y_{j,\beta}$ is also bounded for all $\alpha,\beta = 1,\dotsc,n$.
\end{proposition}

\begin{Proof}
We argue by contradiction.
Indeed, assume that $\fench_{h}(p,y_{j})$ is bounded but $y_{j,\alpha} - y_{j,\beta}\to +\infty$ for some $\alpha,\beta\in\act\equiv\{1,\dotsc,n\}$.
Then, by relabeling indices and passing to a subsequence if necessary, we may assume that
\begin{inparaenum}
[\itshape a\upshape)]
\item
$y_{j,\alpha} \geq y_{j,\kappa} \geq y_{j,\beta}$ for all $\kappa\in\act$;
and
\item
the index set $\act$ can be partitioned into two nonempty sets, $\act^{+}$ and $\act^{-}$, such that $y_{j,\alpha} - y_{j,\kappa}$ is bounded for all $\kappa\in\act^{+}$ while $y_{j,\alpha} - y_{j,\kappa} \to +\infty$ for all $\kappa\in\act^{-}$ (obviously, $\alpha\in\act^{+}$ and $\beta\in\act^{-}$).
\end{inparaenum}
Hence, letting $x_{j} = \choice(y_{j})$, we readily obtain:
\begin{flalign}
\braket{y_{j}}{p - x_{j}}
	&= \sum_{\kappa\in\act} y_{j,\kappa} (p_{\kappa} - x_{j,\kappa})
	= \sum_{\kappa\in\act} (y_{j,\kappa} - y_{j,\alpha}) (p_{\kappa} - x_{j,\kappa})
	\notag\\[1ex]
	&= \sum_{\kappa\in\act^{+}} (y_{j,\kappa} - y_{j,\alpha}) (p_{\kappa} - x_{j,\kappa})
	+ \sum_{\kappa\in\act^{-}} (y_{j,\kappa} - y_{j,\alpha}) (p_{\kappa} - x_{j,\kappa}).
\end{flalign}
The first sum above is bounded by assumption.
As for the second one, Proposition \ref{prop:choice} gives $x_{j,\kappa}\to0$ as $n\to\infty$ for all $\kappa\in\act^{-}$, so $\liminf_{j} (p_{\kappa} - x_{j,\kappa}) > 0$ (recall that $p\in\intsimplex$).
We thus obtain $\sum_{\kappa\in\act^{-}} (y_{j,\kappa} - y_{j,\alpha}) (p_{\kappa} - x_{j,\kappa}) \to -\infty$
and, hence, $\braket{y_{j}}{p - x_{j}} = h(p) - h(x_{j}) - \fench_{h}(p,y_{j}) \to -\infty$, a contradiction.
\end{Proof}

Our final result shows that the Fenchel coupling evolves in a particularly simple fashion under the reinforcement learning dynamics \eqref{eq:RL}:

\begin{lemma}
\label{lem:Fenchel-gradient}
Let $y(t)$ be a solution orbit of \eqref{eq:RL}.
Then, for all $p_{k}\in\strat_{k}$, we have:
\begin{equation}
\label{eq:Fenchel-gradient}
\frac{d}{dt} \fench_{h_k}(p_{k},y_{k}(t))
	= \braket{\payv_{k}(x(t))}{x_{k}(t) - p_{k}}.
\end{equation}
\end{lemma}

\begin{Proof}
By the differential formulation \eqref{eq:RL-diff} of the dynamics \eqref{eq:RL}, we readily obtain:
\begin{flalign}
\frac{d}{dt} \fench_{h_k}(p_{k},y_{k})
	&= \frac{d}{dt} h_{k}^{\ast}(y_{k}) - \braket{\dot y_{k}}{p_{k}}
	\notag\\
	&= \braket{\dot y_{k}}{dh_{k}^{\ast}(y_{k})} - \braket{\payv_{k}}{p_{k}}
	= \braket{\payv_{k}}{\choice_{k}(y_{k}) - p_{k}}
	= \braket{\payv_{k}}{x_{k} - p_{k}},
\end{flalign}
where the first equality in the second line follows from the chain rule and the penultimate one from the fact that $\choice_{k} = dh_{k}^{\ast}$ (Proposition \ref{prop:choice-dh}).
\end{Proof}



\begingroup
\renewcommand{\addcontentsline}[3]{}
\renewcommand{\section}[2]{}

\bibliographystyle{ormsv080}
\bibliography{IEEEabrv,Bibliography-RL}

\endgroup

\end{document}